\documentclass[12pt]{article}

\input xy
\xyoption{all}

\newtheorem{thm}{Theorem}[section]
\newtheorem{prop}[thm]{Proposition}
\newtheorem{lemma}[thm]{Lemma}

\newcounter{ex}[section]

\newcommand{\Example}{\medskip \noindent {\bf
                        Example.} }

\newcommand{\Definition}{\medskip \addtocounter{thm}{1} \noindent {\bf
                         Definition \arabic{section}.\arabic{thm}.} }

\newcommand{\Proof}{\noindent {\sc Proof.\ }}

\newcommand{\Examplenumb}{\medskip \addtocounter{thm}{1}\noindent {\bf
                           Example
                             \arabic{section}.\arabic{thm}.} }

\newenvironment{eq}{\addtocounter{thm}{1}\begin{equation} }{\end{equation}}

\begin{document}
\title{Twists of symmetric bundles}

\author{Ph. Cassou-Nogu\`es, B. Erez, M.J. Taylor }

\date{5 March  2004}

\maketitle

\tableofcontents

\vfill \eject

\noindent {\bf \Large Introduction}
\addcontentsline{toc}{section}{Introduction}
\medskip

\noindent {\it Invariants of general twists.} A symmetric bundle $(E,f)$
over a noe\-the\-rian ${\bf Z}[\frac{1}{2}]$-scheme $Y$ is a vector
bundle $E$ over $Y$ equipped with a symmetric isomorphism $f$
between $E$ and its $Y$-dual $E^{\vee}$. A symmetric bundle can also be
viewed as a quadratic form on $E$ and we  write $(E,q)$
if we take this point of view, or, if the form is clear from the
context, we might even just write $E$. It is well known how to describe
the set of all {\it twists} of $(E,f)$, that is the set of
symmetric bundles which become isomorphic to $(E,f)$ after an \'etale 
base extension. If ${\bf O}(E)$ denotes the orthogonal group
(scheme) attached to $(E,f)$,  this set is $H^1(Y, {\bf O}(E))$
(see [Mi], chapter 3 for a precise definition of this set). For $\alpha$ in $H^1(Y, {\bf O}(E))$,
let $E_{\alpha}$ be the twist of $E$ corresponding to $\alpha$.
Every symmetric bundle of rank $n$ is a twist of the standard
symmetric bundle $(T_n,q_n)=(O_Y^n, x_1^2+...+x_n^2)$. Let ${\bf O}(n)$
denote the automorphism group of this symmetric bundle. We write $\alpha _E$
for the class of $E$ in $H^1(Y, {\bf O}(n))$ ($n=\mbox{rank}(E)$).

Following Delzant \cite{Delz} and Jardine \cite{J-HW}, for any 
symmetric bundle $E$ over $Y$ one can define a cohomological
invariant, which generalizes the classical invariants of quadratic
forms and which is known as the {\it total Hasse-Witt class}. This
is a class $w_t(E)$ in the (graded) \'etale cohomology group
$H^{*}(Y,{\bf Z}/2{\bf Z})$:
$$
w_t(E)= 1+w_1(E)t+w_2(E)t^2+...
$$
A  brief review of the definitions of the
Hasse-Witt invariants can be found in
\cite{CNET1}, Sect.~1.e.
The terms $w_1$ and $w_2$ in degrees one and two
generalize the discriminant and the Hasse-Witt
invariant respectively and have the following elementary description.
 We define  $\delta ^1$  as the map induced by the
determinant map
$$
\delta ^1 = \delta ^{1}_{E}: H^1(Y,{\bf O}(E)) \rightarrow   
H^1(Y,{\bf Z}/2{\bf Z})
$$ and  $\delta ^2$ as the boundary map
$$
\delta ^2 = \delta ^{2}_{E}:   H^1(Y,{\bf O}(E)) \rightarrow   
H^2(Y,{\bf Z}/2{\bf
Z})
$$
associated to the exact sequence of \'etale sheaves of groups
\begin{eq}
1\rightarrow {\bf Z}/2{\bf Z} \rightarrow {\widetilde {\bf O}}(E)
\rightarrow {\bf O}(E) \rightarrow 1 \  , \label{es}
\end{eq} 
\noindent where ${\widetilde {\bf O}}(E)$ is a certain 
covering constructed with the help of the Clifford algebra of $E$ (see  $1$.b for the precise definition). Then
$w_1(E)=\delta ^1_n(\alpha_E)$ and
$w_2(E)=\delta ^2_n(\alpha_E)$,
where,  $\delta ^i_n = \delta ^{i}_{{\bf 1}_n}$.

Let us define a class $\Delta _t(\alpha)$ by the equality
$$
w_t(E_{\alpha}) = w_t(E)\Delta_t(\alpha) \ .
$$
We shall show the  following theorem, which  is a generalisation of a 
result of Serre
\cite{Serre} and  which describes $\Delta _t(\alpha)$ up to terms of
degree $3$.

\begin{thm} \label{general}
Let $(E,q)$ be a symmetric bundle over a scheme $Y$ and let
$\alpha$ be an element of the cohomology set  $H^1(Y,{\bf O}(E))$.
Then
\begin{itemize}
\item[i)] $ w_1(E_{\alpha})=w_1(E)+\delta ^1({\alpha})\ .$

\item[ii)] $ w_2(E_{\alpha})=w_2(E)+w_1(E)\delta ^1({\alpha}) +
\delta ^2({\alpha})\ .$
\end{itemize}
\end{thm}

\noindent Our proof of this result consists of  a cocycle
computation in the spirit of the work of Serre and Fr\"ohlich. It
is based on an explicit formula for the
cup-product of two $1$-cocycles and the study of the behaviour of the
above exact sequence (\ref{es}) under change of forms.
\medskip

\noindent {\it Invariants of Fr\"ohlich twists: the \'etale case.} In the main
body of the paper we will be interested in twists of symmetric
bundles which are obtained from certain types of coverings of $Y$.
Before giving the general definition we describe a special case first 
considered by Fr\"ohlich 
(see \cite{F-SW} and also \cite{EKV}). Let $G$ be a finite group
and let $X$ be a $G$-torsor over $Y$. Consider an orthogonal  representation of $G$ given by a symmetric bundle $(E,q)$ over $Y$ together with a group homomorphism 
$\rho : G\rightarrow \Gamma (Y,{\bf O}(E))$. The $G$-torsor $X$ defines an
element $c(X)$ in $H^1(Y,G)$ and so by push-forward along $\rho$
it defines an element $\rho (X)=\rho_*(c(X))$ in $H^1(Y, {\bf
O}(E))$. We seek a description of $E_{\rho (X)}$.
Let  $Tr_{X/Y}$ denote the bilinear trace form on
$\pi _*(O_X)$.

\begin{prop} \label{torsor-intro}
Consider the $O_Y$-bundle $(E\otimes _{O_Y}\pi _*(O_X))^G$ of
fixed points under $G$, where $G$ acts diagonally: namely through
$\rho$ on the first component and through the given action on the
second. Then
$$
E_{\rho (X)} = E_{\rho ,X} := (E\otimes _{O_Y}\pi _*(O_X))^G
$$
and $q_{\rho ,X}$ may be thought of as the restriction to the
$O_Y$-module $E_{\rho ,X}$ of the form ${\mid G \mid
}^{-1}(q\otimes Tr_{X/Y})$.
\end{prop}

\noindent (See Prop.~\ref{torsor}.)
In comparing the invariants of $E$ and $E_{\rho ,X}$, there 
appear not only the Stiefel-Whitney classes $w_i(\rho)$ of $\rho$ but also   a new
kind of invariant of an orthogonal representation called  the spinor
classes. Building on the work by Fr\"ohlich \cite{F-SW}, Kahn
\cite{K} and Snaith \cite{Snaith}, Jardine showed that,  for $Y$ the
spectrum of a field $K$ of characteristic different from $2$,   there is a class
$$
sp_t(\rho)=1+sp_1(\rho)t +sp_2(\rho)t^2+...
$$
called the {\it total spinor class}, which satisfies
$$
w_t(E_{\rho , X})sp_t(\rho)=w_t(E)w_t(\rho)
$$
in $H^{*}(K,{\bf Z}/2{\bf Z})$ and whose odd components are all
trivial (see \cite{J-HS}). In fact little else is known about the
spinor class except in degree $2$. We extend the work of Serre
\cite{S-HW} and Fr\"ohlich \cite{F-SW} to bundles over schemes. The  proof of our result is different from that of Kahn in  \cite{K} Cor.~6.1.

\begin{thm} \label{etale}
Let $(E,q)$ be a symmetric bundle over a scheme $Y$, let $X$ be a
$G$-torsor over $Y$ and let
$
\rho: G \rightarrow \Gamma (Y, {\bf O}(E))
$
be an orthogonal representation of $G$. Let
$(E_{\rho ,X},q_{\rho ,X})$ be the twist of $(E,q)$ by $\rho$. Then, we 
have the
equalities:
\begin{itemize}
\item[i)] $ w_1(E_{\rho ,X}) = w_1(E)+ w_1(\rho) \ .$
\item[ii)] $ w_2(E_{\rho ,X})=w_2(E)+w_1(E)w_1(\rho)+ w_2(\rho)+ sp_2(\rho)\ .$
\end{itemize}
\end{thm}

\noindent In the proof of this theorem we shall make use of Thm.~\ref{general}.
This gives a new proof of Thm. 2.3  in \cite{EKV} and Thm. 0.2  in \cite{CNET1} in the \'etale case.
Note that our definition
(Def.~\ref{sw2-sp2}) of $sp_2(\rho)$ is different from that in
\cite{K}, but it is known that the two definitions coincide when 
$Y$ is the spectrum of a field, by a remark in \cite{TS1}, p.$127$.
\medskip

\noindent {\it 
Fr\"ohlich twists with tame ramification.} In \cite{S-theta}
Serre considered coverings of Riemann surfaces with odd
ramification and showed how to obtain analogous formul\ae, which
involve expressions defined in terms of ramification data.  Serre's work has
been extended in \cite{EKV}, \cite{K} and more recently in
\cite{CNET1}, \cite{CNET2}. Let
$$
\pi : X\rightarrow Y = X/G
$$
be a covering which is tamely ramified along a divisor $b$ with normal
crossings and whose ramification indices are all odd. In \cite{CNET1} we 
studied the
symmetric bundle $(\pi _*({\cal D}^{-1/2}_{X/Y}),Tr_{X/Y})$, where
${\cal D}^{-1/2}_{X/Y}$ is the locally free sheaf over $X$ whose
square is the inverse of the different of $X/Y$ and $Tr_{X/Y}$ is
the trace form. In this generality we obtained the same formula as
Serre, {\it with no extra terms}. Still the ramification invariant,
which makes its appearance, 
was not so well understood. Here we show that the square root of the
inverse different bundle allows to give an explicit description of
the twists by orthogonal representations coming from tame
coverings and we show how to decompose the ramification invariant
along characters. In our opinion this gives a better understanding of this
invariant.

To be more precise,  let $\pi$ be as above and let $\rho$ be an
orthogonal representation of $G$
$$
\rho :G \rightarrow \Gamma (Y, {\bf O}(E)) \ ,
$$
where $(E,q)$ is a symmetric bundle over $Y$.

\Definition The {\it twist of $(E,q)$ by
$\rho$} is the symmetric bilinear form $(E_{\rho ,X},q_{\rho ,X})$
on $Y$, where
$$
E_{\rho ,X} = (E\otimes _{O_Y}\pi _*({\cal D}^{-1/2}_{X/Y}))^G
$$
is the $G$-fixed submodule of the
$G$-module $E\otimes _{O_Y}\pi _*({\cal D}^{-1/2}_{X/Y})$ and
where $q_{\rho ,X}$ is the form which is the
restriction of ${\mid G \mid }^{-1}(q\otimes Tr_{X/Y})$ to the
$O_Y$-module $E_{\rho , X}$.
\bigskip

\noindent It is clear that this generalizes the \'etale case,
because in that case the inverse different is just $O_X$.

\begin{thm} \label{itsabundle}
\begin{itemize}
\item[i)] The twist $(E_{\rho ,X},q_{\rho ,X})$  of $(E,q)$ by $\rho$ is a
symmetric bundle over $Y$.
\item[ii)] Let $\phi : Z\rightarrow Y$ be a scheme flat over $Y$
and let $T'=Z\times _YX$. For any orthogonal representation $\rho$
of $G$ in ${\bf O}(E)$ we have
$$
(\phi ^{*}(E))_{\rho ,T'} = \phi ^{*}(E_{\rho ,X}) \ .
$$
\end{itemize}
\end{thm}
\bigskip

\noindent We wish to compute the difference
$
w_k(E_{\rho ,X})-w_k(E)
$.
We do not know how to do this directly based on
Thm.~\ref{general} and a cocycle computation as for
Thm.~\ref{etale}, although we suspect it might be possible to do
so. Hence we proceed as in \cite{CNET1} by reducing to the \'etale
case. In order to carry out this reduction we will use a number of functorial
properties of the Hasse-Witt invariants and of the process of twisting,
which are also dealt with in the paper.

Let $G_2$ be a $2$-Sylow subgroup of $G$. We write $Z={ X}/G_2$
and we let  $T$ be  the normalisation of the fiber product
$T'=Z\times _YX$. Hence we have the  diagram.
$$
\begin{array}{ccccl}
T & \longrightarrow & T'= Z\times _YX & \rightarrow & X \\
 \ & & & & \\
 &\pi _Z \searrow & \downarrow & & \downarrow \pi \\
 \ & & & & \\
 &  & Z & \stackrel{\phi}{\rightarrow} & Y
\end{array}
$$
It follows from \cite{EKV}, Sect.~3.4,  that  the action of $G$ on
$X$ induces a $G$-action on $T$ and that $\pi _Z$ identifies $Z$
with $T/G$. Moreover we know from \cite{CNET1}, Thm.~2.2,  that
$\pi _Z: T \rightarrow Z$ is \'etale and hence a $G$-torsor (see [CEPT1], p. 291).
Since the degree of the cover $Z/Y$ is odd, the pull-back map $\phi
^{*}:H ^*(Y, {\bf Z}/2{\bf Z}) \rightarrow H^*(Z, {\bf Z}/2{\bf Z})$
is injective. Therefore we lose no information concerning
the difference $ w_k(E_{\rho , X})-w_k(E)$ by considering the pull-back
 $ \phi ^{*}(w_k(E_{\rho ,X}))-\phi ^{*}(w_k(E))$. We now observe that
 firstly,   we obtain a symmetric bundle over $Z$ by considering  the
 pull-back  $(\phi ^{*}(E), \phi ^{*}(q))$ of $(E,q)$ and secondly,
 that  $\rho $ induces an orthogonal representation $\phi ^{*}(\rho)
 : G \rightarrow \Gamma (Z,{\bf O}(\phi ^{*}(E)))$. 
Hence  we are in a situation where
 we  can use our  construction of a twist and associate to  $(\phi
 ^{*}(E), \phi ^{*}(q))$ the symmetric bundles
 $(\phi ^{*}(E)_{\phi ^{*}(\rho ) , T'},
 \phi ^{*}(q)_{\phi ^{*}(\rho ) , T'})$ and 
$(\phi ^{*}(E)_{\phi ^{*}(\rho ) , T},
 \phi ^{*}(q)_{\phi ^{*}(\rho ) , T})$. Using the
 good functorial properties of the Hasse-Witt invariants and  of the
 process of twisting, the previous difference can be written
 $$
 w_k(\phi ^{*}(E)_{\phi ^{*}(\rho ) , T'}) - w_k(\phi ^{*}(E)) \ .
 $$
Our strategy will be to express this difference as a sum of two terms,
namely:
$$
 \left(w_k(\phi ^{*}(E)_{\phi ^{*}(\rho) , T'})-
w_k(\phi ^{*}(E)_{\phi ^{*}(\rho ) , T})\right) +
\left( w_k(\phi ^{*}(E)_{\phi ^{*}(\rho ) , T})-w_k(\phi
^{*}(E)\right)
  \ .
$$
Since the cover $T/Z$ is \'etale the second term is known by
Thm.~\ref{etale}. Therefore the heart of the problem will be  to compute
the first term. For $w_1$ the formula is as simple as one could 
wish, but for $w_2$
a new class appears. Before stating the result we introduce this new class
$R(\rho, X)$.
\bigskip

\noindent {\it The ramification invariant.}
We start by defining  a divisor on $Y$ which depends on the
decomposition of $\rho$ when restricted to the inertia groups of
the generic points of the irreducible components of the branch
locus $b$ of $\pi$. We number the irreducible components $b_h$
of $b$ by $1\leq h \leq m$, we denote by $\xi_h$ the generic point of $b_h$ and by  $B_h$ 
an irreducible component of the ramification locus of $\pi$ such
that $\pi(B_h)=b_h$. We let $I_h$ be  the inertia group of the
generic point of $B_h$. It follows from our hypotheses that $I_h$
is cyclic of odd order $e_h$. The action of  $I_h$ on the
cotangent space at the generic point  of $B_h$ is given by a
character denoted by $\chi_h$.  For $0 \leq k < e_h/2$ we let
 $ d_k^{(h)}(E)$ denote the rank over the d.v.r. 
$O_{Y,\xi_h} $ of the $\chi_h^{k}$-component of $E_{\xi_h}$ considered 
as an $O_{Y,\xi_h}[I_h]$-module (see section $4$.e for the details).  
We then consider the divisor
$$
R(\rho , X) = \sum _{h=1}^{m}d^{(h)}(E)b_h \ ,
$$
where for $1 \leq h\leq m$ we have put
$$
d^{(h)}(E)=\sum _{k=0}^{e_h/2} kd_k^{(h)}(E) \ .
$$
We shall denote by the same symbol this divisor, the class in
$Pic(Y)$ of the line bundle $O_Y(R(\rho ,X))$ and its image in
$H^2(Y_{et}, {\bf Z}/2{\bf Z})$ under the boundary map
$$
H^1(Y_{et}, {\bf G}_m) \rightarrow H^2(Y_{et}, {\bf Z}/2{\bf Z})
$$
associated to the Kummer sequence
$$
0\rightarrow {\bf Z}/2{\bf Z} \rightarrow {\bf G}_m
\rightarrow {\bf G}_m \rightarrow 0 \
$$
(note that under our assumptions we may  identify $\mu_2$ and
${\bf Z}/2{\bf Z}$).
\medskip

\noindent {\bf Remark.} 1) One recovers the ramification 
invariant in \cite{CNET1} by
considering the regular representation (see Example at the end of Sect.~4).
In \cite{CNET1} we used the notation
$\rho (X/Y)$ for the ramification invariant. This seemed inappropriate here,
as $\rho$ is better used to denote orthogonal representations. (Serre had used
$\omega (X/Y)$, which we changed because it reminded us
too much of  the canonical class...)

2) As pointed out to us by Serre, the ramification
invariant can be viewed as a ``half of a Woods-Hole element''
(our terminology). Namely,
when considering the Lefschetz-Riemann-Roch theorem one encounters
expressions involving terms like $(1-\zeta)^{-1}$, where $\zeta$ is a
root of unity, which come from the inverse of the term
$\lambda_{-1}({\cal N})$ (see for instance \cite{E2} Sect.~2.a, p. 126). 
Now, if $\zeta$ is of order $e$, then
$$
\frac{1}{1-\zeta} = -\frac{1}{e}\sum _{k=1}^{e-1}k\zeta ^k \ .
$$
\bigskip

\noindent {\it Invariants for
Fr\"ohlich twists with tame ramification.} 
We are now in a position to state our next main result.

\begin{thm} \label{tame}
We have the following equalities:
\begin{itemize}
\item[i)] $\ \ \ \ \ w_1(\phi^*(E_{\rho ,X}))= w_1(\phi^*(E))+
w_1(\phi^*(\rho))  \ .$
\item[ii)] $$\begin{array}{rcl}
w_2(\phi^*(E_{\rho ,X})) & = &
w_2(\phi^*(E)) + w_1(\phi^*(E))w_1(\phi^*(\rho)) +
w_2(\phi^*(\rho)) + \\
& & \\
& & + sp_2(\phi^*(\rho)) + \phi^*(R(\rho , X))  \ . 
\end{array}   $$ 
\end{itemize}
\end{thm}
\medskip

\noindent As in \cite{CNET1} the proof of this result relies on two
main ingredients:
on the one hand a formula which expresses the difference between  the total
Hasse-Witt class of the two symmetric bundles
$$
\Upsilon ^{(0)} = (\phi ^{*}(E_{\rho ,X}, q_{\rho ,X}))
\ \ \mbox{and} \ \
\Upsilon ^{(m)} =
(\phi ^{*}(E_{\rho ,T}),(-1)^m\phi ^{*}(q)_{\rho ,T}))^{(-1)^m}
$$
and on the other hand an explicit determinant computation. For reason of simplicity and  when there is no risk of ambiguity we will just denote by $-E$ the symmetric bundle $(E,-q)$. 

As shown in \cite{CNET2} the way to understand the formula for the
Hasse-Witt class is as
a formula expressing the total Hasse-Witt class of
a {\it metabolic symmetric complex} in terms of Chern classes of a
lagrangian, that is a 
maximal totally isotropic subcomplex of this complex
(see below for the terminology). More
precisely, by decomposing the normalisation map
$T\rightarrow T'$  into a sequence of normalisations, 
where we add in one irreducible component at a time, we show how 
to construct
a lagrangian of $\Upsilon ^{(0)}\perp -\Upsilon ^{(m)}$ viewed
as a symmetric complex concentrated in degree $0$. 
We consider the sequence of $Z$-morphisms
$$
T=T^{(m)}\rightarrow T^{(m-1)}\rightarrow \cdots \rightarrow T^{(0)}=T' \ ,
$$
numbered by the $m$ components of the branch locus $b$ of
the covering and obtained
by normalisation along a component of $b$ 
as described in \cite{CNET1}, Sect.~3.
We will prove in Sect.~4 that,  for $1 \leq h \leq m$,
we obtain an exact sequence of locally free $O_Z$-modules
$$
0\rightarrow (\phi ^{*}(E)\otimes _{O_Z} I^{(h)})^G\rightarrow
\phi ^{*}(E)_{\rho , T^{(h)}}
\oplus   \phi ^{*}(E)_{\rho , T^{(h+1)}}
\rightarrow (\phi ^{*}(E)\otimes _{O_Z}{\cal G}^{(h)})^G \rightarrow
0 
$$
and that  for $0\leq h\leq m-1$  the symmetric bundle defined as the
orthogonal sum $(\phi ^{*}(E)_{\rho , T^{(h)}},
\phi ^{*}(q)_{\rho , T^{(h)}})\perp
(\phi ^{*}(E)_{\rho , T^{(h+1)}},-\phi ^{*}(q)_{\rho , T^{(h+1)}})$
is metabolic.  

For a bundle $V$ of rank $n$ over Y we denote by $c_i(V)$ the $i$-th Chern
class of $V$ as an element of
$H^{2i}(Y,{\bf Z}/2{\bf Z})$, [Gr].  We define an element of
$H^*(Y, {\bf Z}/2{\bf Z})$ by
$$
d_t(V) = \sum _{i=0}^{n}(1 +(-1)t)^{n-i}c_i(V)t^{2i} \ .
$$
To understand fully the next result it is necessary to work in the derived category of bounded complexes of vector bundles over a scheme. 
The reader is referred to \cite{Ba1}, \cite{Ba2} and \cite{Wa}, 
for the basic theory and 
to [CNET2] for a discussion in our 
geometric context (see also the Appendix). 
The key point is that while $\Upsilon ^{(0)}\perp
-\Upsilon ^{(m)}$ is not a metabolic symmetric bundle, it is,
however, a metabolic complex in the derived category, if we view
it as a symmetric complex
concentrated in degree $0$. Furthermore it then has for  lagrangian
the complex  $M_{\bullet}$ which is concentrated in degrees $1$ and $0$
with $M_0= \oplus_{h=0}^{h=m-1} {(\phi^*(E)\otimes I^{(h)})}^G$ and
$ M_1= \oplus_{h=1}^{h=m}\phi^*(E)_{\rho,T^{(h)}} $.    

Using the generalisation of the main lemma of Sect.~4.b
to derived categories we
then  get 
$$ w_t(\Upsilon ^{(0)})w_t(-\Upsilon ^{(m)}) = d_t(M_0 - M_1) \ .
$$
This formula lies behind the following result. 

\begin{thm} \label{big-main-lemma}
In $H^*(Z,{\bf Z}/2{\bf Z})$ the class 
$w_t(\phi ^{*}(E_{\rho , X}, q_{\rho , X}))$
equals
$$
w_t(\phi ^{*}(E)_{\rho , T},(-1)^m\phi ^{*}(q)_{\rho , T})^{(-1)^m}
\prod _{0 \leq h \leq m-1} d_t((\phi ^{*}(E)
\otimes _{O_Z}
{\cal G}^{(h)})^G)^{(-1)^h} \ .
$$
 \end{thm}

\noindent Theorem 0.7  results from this by making
the degree two terms explicit. It is quite  remarkable to see
how the ramification invariant comes out from this.
\medskip

\noindent To end this introduction let us point out that: firstly
Theorem~\ref{tame}
is new even for curves; secondly it provides a substantial strengthening 
of the main result of [CNET1], since we no longer need to 
impose such stringent regularity conditions. Indeed, 
fix  a subgroup $H$ of $G$ and let
$$ 
\lambda: X \rightarrow V:=X/H
$$ 
denote the quotient map and 
$$
\gamma: V \rightarrow Y
$$ 
the induced map, so that $\pi= \gamma \circ \lambda$. 
Note that $V$, being the  quotient of a normal scheme by a finite group, 
is normal but not necessarily regular.
Also observe that when $\lambda$ is flat,   then
the symmetric bundle $(\gamma_*({\cal D}^{-1/2}_{V/Y}), Tr_{V/Y})$ which, 
under this assumption, was  the main object of investigation of \newline
 [CNET1], 
provides  an example of a twist of a symmetric bundle by an orthogonal 
representation which is  both natural and explicit (see Example at the end of
Sect.~4). So our previous work  does become
 a particular case of the general study of this paper. 
Moreover, when $V$ is not regular and hence  $\lambda$ not flat, then   
$(\gamma_*({\cal D}^{-1/2}_{V/Y}), Tr_{V/Y})$ is not in general a 
symmetric bundle. In this case   
$(\pi _*({\cal D}^{-1/2}_{X/Y})^H,Tr_{V/Y})$  provides the right substitute.

We should also indicate some  other related developments.
The twisting results for forms over fields  play
an important role in the work of Saito on the
sign of the functional equation of the L-function of an orthogonal
motive (see \cite{TS1}). Saito proves a $p$-adic version of
Fr\"ohlich's result: namely,  he deals with Galois representations
that do not necessarily have finite image.  He uses this 
to prove a result analogous to the Fr\"ohlich-Queyrut Theorem, which
states  that the global root number of  real  orthogonal
characters equals one \cite{FQ}. Saito's approach follows that of
Deligne, who interpreted this
reciprocity result in terms of the Stiefel-Whitney classes
of the (local) characters \cite{De}.
More recently
the results of \cite{CNET1} have been
used by Glass in \cite{Glass} to relate the Galois invariants appearing
there to $\epsilon$-factors. On a different track, Saito has generalized
Serre's original formula to the case of (smooth) non-finite morphisms
$X\rightarrow Y$ \cite{TS2}.
The conjunction of these  results
suggests a beautiful picture, parts of which  are still hidden,
in which the formul{\ae}  obtained in this paper seem to
hold a special place.
\bigskip

\section{General twists}

\subsection{Symmetric bundles on schemes} \label{Quadratic}
For completeness, we recall the basic definitions concerning forms
over sche\-mes.
We let $Y$ be a scheme and we assume that $2$ is invertible over
$Y$, then the theory of symmetric bilinear forms over $Y$ is equivalent to
that of quadratic forms over $Y$.  A {\it vector bundle} $E$ on $Y$ is a
locally free $O_Y$-module of finite rank. The dual of a vector
bundle $E$ is the vector bundle $E^{\vee}$ such that, for any open
subscheme $Z$ of $Y$
$$
E^{\vee} (Z) = Hom _{O_Z}(E|_Z, O_Z) \ .
$$
There is a natural evaluation pairing $<\ , \ >$ between $E$ and
$E^{\vee}$ and one can identify $E$ with the double dual $E^{\vee
\vee}$ by
$$
\kappa : E \cong E^{\vee \vee} \ ,
$$
where $<\alpha, \kappa (u)> = <u, \alpha >$. A {\it symmetric
bilinear form}  on $Y$ is a vector bundle $E$ on $Y$ equipped with
a map of sheaves
$$
q : E \times _Y E\rightarrow O_Y \ ,
$$
which on sections over an open subscheme
 restricts to a symmetric bilinear form.
This defines an {\it adjoint} map
$$
\varphi = \varphi _q : E\rightarrow E^{\vee} \ ,
$$
which because of the symmetry assumption equals its transpose:
$$
\varphi = \varphi ^t : E\stackrel{\kappa}{\rightarrow} E^{\vee
\vee} \stackrel{\varphi ^{\vee}}{\rightarrow} E^{\vee} \ .
$$
We shall say that $(E,q)$ is {\it non-degenerate} (or unimodular)
if the adjoint $\varphi$ is an isomorphism. From now on we will
call a {\it symmetric bundle} any vector bundle endowed with a
non-degenerate quadratic form.

\subsection{Clifford algebras}

\noindent The properties of the Clifford algebra and the Clifford group
associated with a symmetric bundle will play an important role in
the proof of Thm.~\ref{general}. Therefore we start  by  briefly
recalling some of  the basic material that we shall need in this section. Our
references will be \cite{EKV}, section 1.9, [Knu], Chapter 4, in the
case  of  forms over a ring, and [F], Appendix 1, for a  brief
review in the case of forms over a field. In fact one has to
observe that  most  of the definitions about Clifford
algebras associated to a quadratic form over a field,  or more
generally over a commutative ring,  can be generalised in our
geometric context. Moreover,   by reducing  to  affine
neighbourhoods,   we will essentially work with non-degenerate
forms over  noetherian,  integral domains.

To any symmetric bundle $(V,q)$ of constant rank $n$,
 one associates a sheaf of
algebras $\cal C$$ (q)$ over $O_Y$,  of constant rank $2^n$.
 As in the classical case one
has the notion of odd and even elements of $\cal C$$ (q)$ and hence
a  $\bmod  2$ grading. The Clifford group ${\cal C}^*$$(q)$ is the
subgroup of homogeneous, invertible elements $x$ in  $\cal C$$(q)$
such that $xvx^{-1}$ belongs to $V$ for any $v$ of $V$. The
$\bmod 2$ grading induces a splitting
$$ 
{\cal C}^*(q)={\cal C}^*_{+}(q)\cup {\cal C}^*_{-}(q) \ .
$$
Let $\sigma $ be the anti-automorphism on ${\cal C}(q)$ induced by the
identity on $V$ so that $\sigma (v_1\cdots v_m)=v_m\cdots v_1$. One verifies
that the map  $N$ defined on $ {\cal C}^*(q)$ by $N(x)=\sigma(x)x$
induces an homomorphism,
$$
N: {\cal C}^*(q) \rightarrow {\bf G}_m \ .
$$
We define the algebraic group scheme $ \widetilde {{\bf O}}(q)$ as
the kernel of this homomorphism. This group scheme also
splits as  $ \widetilde {{\bf O}}_{+}(q)\cup  \widetilde {{\bf
O}}_{-}(q)$. Let $x$ be in  
$ \widetilde {{\bf O}}_{\epsilon}(q)$ with $\epsilon =\pm 1$, 
then we define $r_q(x)$
as the element of ${\bf O}(q)$
$$
\begin{array}{rcl}
r_q(x) :  V & \rightarrow & V \\
v & \mapsto & \epsilon xvx^{-1} \ .
\end{array}
$$
This defines a group homomorphism $r_q:\widetilde {{\bf O}}(q)
\rightarrow {\bf O}(q)$. One can show  that for each $x \in
\widetilde {{\bf O}}_{\epsilon}(q)$ the element  $r_q(x)$ belongs
to $ {\bf O}_+(q)={\bf SO}(q)$ or ${\bf  O}_-(q)={\bf O}(q)
\setminus  {\bf SO}(q)$ depending on whether  $\epsilon=1$ or
$-1$. Where there is no risk of ambiguity, we will write $r$ for
$r_q$.

We then have constructed  an exact sequence of \'etale sheaves of
groups:
$$ 1\rightarrow {\bf Z}/2{\bf Z} \rightarrow {\widetilde {\bf O}}(q)
\rightarrow {\bf O}(q) \rightarrow 1 \ .$$ We recall that we have
previously introduced
$$\delta ^1: H^1(Y,{\bf O}(q)) \rightarrow   H^1(Y,{\bf Z}/2{\bf Z})$$
as  the map induced by the determinant and
$$\delta ^2:   H^1(Y,{\bf O}(q)) \rightarrow   H^2(Y,{\bf Z}/2{\bf Z})$$
as the boundary map associated to the above exact sequence.
\smallskip
We next consider an affine situation, namely $Y=\mbox{Spec} (R)$ where $R$
is an integral domain. For a symmetric bundle $(V,q)$ over $O_Y$,
by abuse of notation,  we shall write $V$, ${\bf  O}(q)$,
${\widetilde {\bf O}}(q)$, for the corresponding  module  
or groups of the global
sections of these objects. For any invertible element $a$ of
${\cal C}(q)$ we will denote by $\imath _a$ the inner automorphism
of  ${\cal C}(q)$ given by conjugation, $x \mapsto axa^{-1}$.

\begin{prop} \label{clifford} 
Let $(E,q)$ and $(F,f)$ be symmetric 
bundles over $O_Y$ and let $\theta : (E,q) \rightarrow (F,f)$ be an 
isometry. Then
\begin{itemize}
\item[i)] $\theta$ extends to  an isomorphism $\tilde \theta$ from 
$\widetilde { {\bf O}}(q)$ onto $\widetilde {{\bf O}}(f)$ which 
induces the isomorphism $u \mapsto \theta u 
{\theta}^{-1}$ from $\mbox{\rm Im}(r_q)$ onto $\mbox{\rm Im}(r_f)$.

\item[ii)] Suppose that $(E,q)=(F,f)$ and that $\theta $ belongs to 
$\mbox{\rm Im}(r_q)$. 
Let $t(\theta)$ denote a lift of $\theta$ to  $\widetilde {{\bf O}}(q)$.
If $\theta$ belongs to ${\bf O}_+(q)$ (resp. ${\bf O}_-(q)$), then
on $\widetilde {{\bf O}}_{\epsilon}(q)$ there holds
$\tilde \theta =\imath_{t(\theta)}$ (resp. 
$\tilde \theta = \epsilon \imath_{t(\theta)}$). 
\end{itemize}
\end{prop}

\Proof The universal property of the
Clifford algebra implies  that $\theta$ induces a graded isomorphism
$\tilde\theta: {\cal C} (q) \rightarrow {\cal C} (f)$, which by
restriction induces the required  isomorphism. Moreover since it
coincides with the identity on ${\bf Z}/2{\bf Z}$, it follows that  
$\tilde \theta$ induces a group isomorphism
\noindent $s_\theta: \mbox{Im}(r_q) \rightarrow \mbox{\rm Im}(r_f)$. 
Let $u \in \mbox{\rm Im}(r_q)$,
with $u \in {\bf O}_{\epsilon}(q)$ and let $a  \in  \widetilde
{{\bf O}}_{\epsilon}(q)$ denote a lift of $u$. From the very
definition of $r$ we deduce that for any $y \in F$
$$
s_{\theta}(u)(y)=\epsilon \tilde \theta (a) y \tilde \theta (a^{-1}) \ .
$$
Since $ \tilde\theta$ is an isomorphism of $O_Y$-algebras, the right
hand side of this equality can be written
$$
\epsilon \tilde \theta(a {\tilde\theta}^{-1}(y)a^{-1})=
\theta (\epsilon a\theta^{-1}(y)a^{-1})=(\theta u \theta^{-1})(y) \ .
$$
Hence we have proved that $s_\theta$ and $u \mapsto \theta u\theta ^{-1}$
coincide on $\mbox{Im}(r_q)$.

Under the hypothesis of $(ii)$  we obtain  two automorphisms of
${\cal C}(q)$, namely $\tilde \theta $ and $\imath _{t(\theta)}$.
Moreover,  since $t(\theta) \in \widetilde {{\bf
O}}_{\epsilon}(q)$, it follows from the definition of $r_q$  that
$\theta (x)=\epsilon t(\theta)xt(\theta)^{-1}$ for any $x$ in $E$.
If $\theta \in {\bf O}_+(q)$, then $\epsilon=1$ and $\tilde \theta $
and $\imath _{t(\theta)}$,   which coincide on $E$, coincide on
$\cal C$$(q)$. If   $\theta \in {\bf O}_-(q)$, then $ \epsilon=-1$.
Therefore  $\tilde \theta$ and $\imath _{t(\theta)}$ will coincide
on $\cal C_+$$(q)$ and will differ by a minus sign on  $\cal
C_-$$(q)$, so the result follows.
\bigskip

\subsection{Proof of Theorem~\ref{general}}

\noindent Let $(E,q)$
be a symmetric bundle and let  $\alpha$ be a  class of the cohomology set
$ H^1(Y,{\bf O}(q))$. Since the set
$H^1(Y,{\bf O}(q))$ classifies the isometry classes of twisted
forms of  $(E,q)$, we may consider  a twisted form
$(E_{\alpha}, q_{\alpha})$ whose class represents  $\alpha$.
Finally we denote by $(T_n,q_n)$ the standard/sum of squares form $(O_Y^n,
x_1^2+ \cdots +x_n^2)$ and as usual we write ${\bf O}(n)$ for
${\bf O}(q_n)$. Since both $(E,q)$ and $(E_{\alpha},q_{\alpha})$
are twisted forms of $(T_n,q_n)$, there exists an affine covering
$\cal U =(U$$_i \rightarrow Y)$ for the \'etale topology and
isometries
$$
\varphi_i: (E_{\alpha},q_{\alpha})\times U_i \rightarrow (E,q)\times U_i
$$
$$\psi_i:(E,q)\times U_i \rightarrow (T_n,q_n)\times U_i \ .
$$
Therefore, following [Mi], Sect.~4, we deduce that
$(\alpha_{ij})=\varphi_i{\varphi_j}^{-1}$ and
$(\gamma_{ij})=\psi_i{\psi_j}^{-1}$ are  $1$-cocycles
representing $(E_{\alpha},q_{\alpha})$ in
$H^1({\cal U}/Y,{\bf O}(q))$ and $(E,q)$ in
$H^1({\cal U}/Y, {\bf O}(n))$ respectively. By
considering
$(\delta_{ij})=(\psi_i\varphi_i){(\psi_j\varphi_j)}^{-1}$ we
obtain a $1$-cocycle representative of $(E_{\alpha},q_{\alpha})$
in $H^1(\cal U$$/Y, {\bf O}(n))$. We
observe that we can write
$$
\delta_{ij}=\psi_i\varphi_i\varphi_{j}^{-1}\psi_{j}^{-1}=
\psi_i\alpha_{ij}\psi_{j}^{-1}
=
\psi_i\psi_{j}^{-1}\psi_j\alpha_{ij}\psi_{j}^{-1}=
\gamma_{ij}(\psi_j\alpha_{ij}\psi_{j}^{-1}) \ .
$$
 In order to obtain a representative of $w_1(E_{\alpha})$, it suffices 
to take the image by the determinant map  of the cocycle  $(\delta_{ij})$. 
>From the previous equalities  we deduce that 
$det(\delta_{ij})=det(\gamma_{ij})det(\alpha_{ij})$ which immediately 
implies that
 $$w_1(E_{\alpha})=w_1(E)+\delta^1(\alpha) \ , $$
as required.

We now want to compare $w_2(E_{\alpha})$ and $w_2(E)$.   After
refining $(U_i)$ we may assume, [Mi], III.2.19, that each $
\alpha_{ij}$, (resp.  $\gamma_{ij}$), is the image of an element
of $\widetilde {{\bf O}}(q)(U_{ij})$, (resp. $\widetilde {{\bf
O}}(n)(U_{ij}))$ that we denote by $ \alpha_{ij}'$, (resp.
$\gamma_{ij}'$).  By abuse of notation for any $l \in \{i,j,k\}$
we will still denote by $\psi_l$ the restriction to $U_{ijk}$ of 
the isometry $\psi_l$. We now deduce from Pro. 3.2. $(i)$
that each $\delta_{ij}$ is the image of the element
$\delta_{ij}'=\gamma_{ij}'\tilde  \psi_j (\alpha_{ij}')$. Therefore
$w_2(E_{\alpha})$ is the class of the $2$-cocycle $(b_{ijk})$
where
$$
b_{ijk}= \delta_{jk}'{\delta_{ik}'}^{-1} \delta_{ij}' \in 
\Gamma( U_{ijk}, {\bf Z}/2{\bf Z}) \ .
$$
Using repeatedly the previous equalities we obtain that
 $$
b_{ijk}= \gamma_{jk}'\tilde \psi_k (\alpha_{jk}')\tilde \psi_k 
({\alpha_{ik}'}^{-1})
{\gamma_{ik}'}^{-1}\gamma_{ij}'\tilde \psi_j (\alpha_{ij}') \ ,
$$
that we write
$$
b_{ijk}= \gamma_{jk}'\tilde \psi_k (\alpha_{jk}')\tilde \psi_k 
({\alpha_{ik}'}^{-1})
{\gamma_{jk}'}^{-1}(\gamma_{jk}'{\gamma_{ik}'}^{-1}\gamma_{ij}')\tilde
\psi_j (\alpha_{ij}') \ .
$$ 
We now observe that on the one hand 
$(\gamma_{jk}'{\gamma_{ik}'}^{-1}\gamma_{ij}')=\pm 1$ and thus this factor
commutes with every factor of the product, while on the other hand,
 since $(\alpha_{jk}'{\alpha_{ik}'}^{-1}\alpha_{ij}')=\pm 1$,
we have  $\tilde \psi_k (\alpha_{jk}') \tilde \psi_k
({\alpha_{ik}'}^{-1})=(\alpha_{jk}'{\alpha_{ik}'}^{-1}\alpha_{ij}')\tilde
\psi_k ({\alpha_{ij}'}^{-1})$. Piecing these observations
together we obtain the equality:
$$
b_{ijk}=(\gamma_{jk}'{\gamma_{ik}'}^{-1}\gamma_{ij}')
(\alpha_{jk}'{\alpha_{ik}'}^{-1}\alpha_{ij}')\gamma_{jk}'\tilde \psi_k 
({\alpha_{ij}'}^{-1}){\gamma_{jk}'}^{-1}
\tilde \psi_j (\alpha_{ij}') \ .
$$ 
We consider this equality as a
product of three factors namely
$$ 
b_{ijk}=(\gamma_{jk}'{\gamma_{ik}'}^{-1}\gamma_{ij}')
(\alpha_{jk}'{\alpha_{ik}'}^{-1}\alpha_{ij}')\epsilon_{ijk}
$$
 where $\epsilon_{ijk}=\gamma_{jk}'\tilde \psi_k 
({\alpha_{ij}'}^{-1}){\gamma_{jk}'}^{-1}\tilde \psi_j (\alpha_{ij}')\ .$
The first two factors are  $2$-cocycles which respectively
represent $w_2(E)$ and $\delta^2(\alpha)$. 

We now want to simplify
the expression for  $\epsilon_{ijk}$. 
This is achieved by considering the various possible  
signs of  $det( \gamma_{jk})$ and $det(\alpha_{ij})$. 
 We first observe that the
equality $\gamma_{jk} \psi_k=\psi_j$ implies that $\tilde
\gamma_{jk}\tilde \psi_k=\tilde \psi_j$. Moreover,  with the notation
of Prop.~\ref{clifford},  we may write $\gamma_{jk}'\tilde \psi_k
(x){\gamma_{jk}'}^{-1}=\imath_{\gamma_{jk}'}\tilde \psi_k (x)$. We start by  
supposing  that $det(\gamma_{jk})=1$. Then it follows  Prop.~\ref{clifford} $ii$) 
that $\imath_{\gamma_{jk}'}=\tilde \gamma_{jk}$, hence
$\gamma_{jk}'\tilde \psi_k (x){\gamma_{jk}'}^{-1}= \tilde \gamma_{jk}
\tilde \psi_k (x)=\tilde \psi_j(x)$ and we conclude that
$\epsilon_{ijk}=1$ . We now suppose  that
$det(\gamma_{jk})=-1$.  First we assume that $det(\alpha_{ij})=1$. Then, in this case,   $\alpha_{ij}' \in \cal C_+ $$(q)$ and therefore $\tilde \psi_k ( \alpha_{ij}'^{-1})\in \cal C_+$$(q_n)$.  Using as before   Prop.~\ref{clifford} $ii$),  we obtain that 
$$\gamma_{jk}'\tilde \psi_k (\alpha_{ij}'^{-1}){\gamma_{jk}'}^{-1}= \tilde \gamma_{jk}
\tilde \psi_k (\alpha_{ij}'^{-1})=\tilde \psi_j(\alpha_{ij}'^{-1})$$
 and we conclude that $\epsilon_{ijk}=1$. We now suppose that  $det(\alpha_{ij})=-1$ then $\tilde \psi_k ( \alpha_{ij}'^{-1})\in \cal C_-$$(q_n)$. Therefore, using once again  Prop.~\ref{clifford} $ii$), we deduce that 
$$\gamma_{jk}'\tilde \psi_k (\alpha_{ij}'^{-1}){\gamma_{jk}'}^{-1}= -\tilde \gamma_{jk}
\tilde \psi_k (\alpha_{ij}'^{-1})=-\tilde \psi_j(\alpha_{ij}'^{-1}) \ .$$ It then follows from the definition that $\epsilon_{ijk}=-1 $ .  As an immediate consequence of  the study of these different cases  we obtain  the equality
$$ 
\epsilon_{ijk}=(-1)^{\epsilon (det(\alpha_{ij}))\epsilon (det(\gamma_{jk}))}\ 
.
$$ 
where for $x$ in $\{ \pm 1\}$ we define $\epsilon (x) \in {\bf Z}/2{\bf Z}$ by 
$x=(-1)^{\epsilon (x)}$. Therefore we conclude that $(\epsilon_{ijk})$ is a 
$2$-cocycle representative of the cup product $ w_1(E)\delta^1(\alpha)$. 
This completes the proof of the theorem.
\bigskip

\section{The twist of a bundle {\em \`a la} Fr\"ohlich}

\noindent The aim of this section is to prove Theorem~\ref{itsabundle}: namely, that 
the twist of symmetric bundle by an orthogonal representation
of the ``Galois group'' of a tame covering
is again  a symmetric bundle. It should be observed that in this section we
are able to relax some of the hypotheses imposed in subsequent sections.

\subsection{Tame coverings with odd ramification} \label{Tame}

\noindent In what follows all the schemes will be defined over
$\mbox{Spec} ({\bf Z}[\frac{1}{2}])$. We let $ X$ be a connected,
projective, regular scheme which is  either defined over the
spectrum $\mbox{Spec} ( {\bf F}_p)$ of the prime field of
characteristic $p\not= 2$ or is flat over $\mbox{Spec} ({\bf
Z}[\frac{1}{2}])$. We assume  that $ X$ is equipped with a tame
action by a finite group  $G$, in the sense of Grothendieck-Murre. In particular
the quotient
$$
{ \pi} : { X}\rightarrow Y = { X}/G
$$
exists and  is a torsor for  $G$ over $Y$ outside a divisor $b$
with normal crossings (see [Gr-M], [CEPT1],
 [CE] Appendix and [CEPT2] 1.2 and Appendix). We assume that $\pi$ is
a flat morphism of schemes and therefore  that $Y$ is regular. The
different  ${\cal D}_{{ X}/Y}$ is defined as the annihilator of
the sheaf  $\Omega ^1_{{ X}/Y}$ of relative differentials of
degree $1$(see [Mi]  Rem.~I.3.7). The reduced closed subscheme of
${ X}$ defined by the support of the  different is the ramification locus of $
\pi$, which we denote by ${ B} = B({ X}/Y)$. Then $b = b({ X}/Y)$
is  the reduced sub-scheme of $Y$ defined by the image of ${ B}$
under  $ \pi$. There are decompositions
$$
b({ X}/Y) = \coprod _{1\leq h \leq m}b_h \ \ \ \ \mbox{and}\ \ \ \
B({ X}/Y) = \coprod _{h,k} { B_{h,k}} \ ,
$$
where the $b_h$ are the irreducible components of  $b$ and where
for any fixed integer $h$ between $1$ and $m$, the ${ B_{h,k}}$
are the irreducible components of $ B$ such that ${ \pi}({
B_{h,k}}) = b_h$. The tameness assumption on the ramification
implies that the branch locus $b({ X}/Y)$  on ${ Y}$ is a
divisor with  normal crossings. For each irreducible component
$b_h$ of $b$, we denote by  $\xi _h$ the generic point of $b_h$
and by  ${ \xi}_{h,k}^{'}$  a generic point of the component ${
B_{h,k}}$ of the ramificaton locus on ${ X}$. The ramification index
(resp. residue class extension degree) of $\xi _{h,k}^{'}$ over
$\xi _h$ which is independent of $k$  will be denoted by $e(\xi
_{h}^{'})$ (resp. $f(\xi _{h}^{'}))$. We  assume  that the inertia
group of each closed point of ${ X}$  has odd order, thus  every
point has odd inertia and in particular the integers $e(\xi
_{h}^{'})$ are odd. We introduce the  divisor of $X$
$$\omega_{X/Y}=\sum_{h,k}(e(\xi _{h}^{'})-1) B_{h,k} \ $$
and we  define the square root of the codifferent as a  vector
bundle on $X$ by setting
$$  {\cal D}_{X/Y}^{-1/2}=O_X(\omega_{X/Y}/2)\  .$$
 Therefore we obtain a symmetric bundle on $Y$ as defined in 1.a by considering
$(\pi_*({\cal D}_{ X/Y}^{-1/2}),Tr_{X/Y})$, where $ Tr_{X/Y}$
denotes the trace form.

\subsection{Definition of the twist: the bundle}

\noindent Let $\rho$ be an orthogonal representation
of $G$ taking values values in 
the orthogonal group ${\bf O}(E)$ of the symmetric bundle
$E$ over $Y$. Under the
assumptions that the action of $G$ is tame and that the
ramification indices are odd,  we may consider the $O_Y$-module
$$
E_{\rho , X} := (E\otimes _{O_Y}\pi _*({\cal D}^{-1/2}_{X/Y}))^G \ .
$$
Our  first aim is to show  that $E_{\rho , X}$ is a vector bundle
over $O_Y$. The proof  will be a consequence of 
a number of elementary algebraic
results.
Later on,  after giving a precise definition of the
form $q_{\rho , X}$,  we shall show  that $(E_{\rho , X}, q_{\rho , X})$
is non-degenerate and, in the \'etale case,
 we identify this form to a twist of $E$ by a cocycle (see Prop.~\ref{torsor}).
\bigskip

\subsection{Algebraic preliminaries}

 \noindent  For these auxiliary steps,  as in the rest of this paper,  we adopt the following conventions:
 $R$ is  an integral domain and $M$ is a left  $R[G]$-module
which  we assume to be locally free and finitely generated as an
$R$-module. We write $\sigma=\sum_{g\in G}g$.

\begin{lemma} \label{fixedtrace}
If $M$ is a projective $R[G]$-module, then $M^G=\sigma M$.
Furthermore:
\begin{itemize}
\item[i)]$ M^G$ is a locally free $R$-module.
\item[ii)] The map $m \mapsto \sigma m $ induces an isomorphism of $R$-modules  from $M_G$ onto $M^G$.
\end{itemize}
\end{lemma}

\Proof All the above statements are clear when $M=R[G]$; they
therefore hold for any finitely generated free  $R[G]$-module and
are then easily seen to hold for a direct summand of a finitely
generated free  $R[G]$-module.
\bigskip

\noindent We now consider a symmetric bundle $(M,t)$ over $R$; that is to
say $M$ is a finitely generated locally free $R$-module,  equipped
with a non-degenerate symmetric bilinear form $t$. We denote by
$\varphi_t$ the adjoint map $ \varphi_t:M \rightarrow M^{\vee}$.
We suppose further that $M$ is now  a projective $R[G]$-module and
that the pairing $t$ is $G$-invariant. Under these assumptions  we
then can use Lemma~\ref{fixedtrace}  to define the symmetric bilinear form
$t^G$ on $M^G$ by
$$
t^G(x,y)=t(m,y)
$$
where $m$ is any arbitrary element of $M$ such that $\sigma m=x$.
Let $I$ be the $R$-submodule of $M$ generated by the set $\{(1-g)m,
m \in M,g \in G\}$.  Since $m$ is defined up to an  element of
$I$ (Lemma~\ref{fixedtrace} (ii))  and since $t$ is $G$-invariant, one verifies
immediately that $t^G$ is well defined. Moreover one observes that
for any $x$ and $y$ in $M^G$ one has
$$\mid G\mid t^G(x,y)=t(x,y).  $$

\begin{prop} \label{CohAl}
 If  $M$ is a projective $R[G]$-module, then $(M^G,t^G)$ is a symmetric
 $R$-bundle.
\end{prop}

\Proof From Lemma~\ref{fixedtrace} we know that $M^G$ is a locally free $
R$-module. It remains to prove that $\varphi_{t^G}$ is an
$R$-module isomorphism from $M^G$ onto $Hom_R(M^G,R)$.
>From the exact sequence
$$
0 \rightarrow  I  \rightarrow  M  \rightarrow M_G \rightarrow  0 \ ,
$$
we deduce the exact sequence
$$ 0 \rightarrow  Hom_R(M_G,R)  \rightarrow Hom_R( M,R)
\rightarrow Hom_R(I,R) \rightarrow  Ext_R(M_G,R) \ .  $$
By  Lemma~\ref{fixedtrace} we know that $M_G$ is isomorphic to $M^G$ and is
therefore $R$-projective. Hence we conclude that $Ext_R(M_G,R)=\{0\}$.
We consider the following diagram
$$
\begin{array}{ccccccccc}
0& \rightarrow & M^G & \rightarrow & M & \rightarrow &N &\rightarrow & 0\\
& &\varphi'  \downarrow \ \ \ \ & & \downarrow  \varphi  & &\ \ \
\downarrow \varphi''
& &  \\
0& \rightarrow & Hom_R(M_G,R)  &

\rightarrow & Hom_R(M,R) & \rightarrow & Hom_R(I,R) &\rightarrow
&0
\end{array}
$$
where $N=M/M^G$, $\varphi=\varphi_t$ and $\varphi'$ is
defined for
$x \in M^G$ and $y \in M_G$ by
$$
\varphi'(x)(y)=\varphi(x)(m) \ ,
$$ where  $m$ denotes  any
representative of $y$ in $M$. Once again, since $m$ is defined
up to an element of $I$, $\varphi'$ is well-defined and moreover
the first square and thus the diagram itself are both commutative.
Applying the snake lemma to the previous diagram we obtain the exact sequence
$$
0 \rightarrow  \ker\varphi'  \rightarrow \ker\varphi
\rightarrow \ker\varphi'' \rightarrow  \mbox{\rm coker}\varphi'
\rightarrow \mbox{\rm coker}\varphi \ .
$$ Since $\varphi$ is an
isomorphism, the previous sequence reduces to
$$
0 \rightarrow  \ker\varphi''  \rightarrow  \mbox{\rm coker}\varphi'  
\rightarrow  0 \ .
$$
Since $M_G$ is isomorphic to $M^G$ and $\varphi'$ is injective,
$\varphi'(M^G)$ is an $R$-submodule of $Hom_R(M_G,R)$ of the same
rank. Hence $\mbox{\rm coker}\varphi'$ is an $R$-torsion module. We now  show
that $N$ and hence $\ker\varphi''$ is torsion free. Once again we can reduce consideration
to the case where $M$ is a free $R[G]$-module with
$\{e_1,e_2,...,e_n\}$ as a basis. Let $m=\sum_{1\leq i\leq n} a_ie_i$
 be an element of $M$, and let $d\in R\setminus\{0\}$ be  such
that $dm \in M^G$. Since $ \{\sigma e_1, \sigma e_2,...,\sigma e_n\}$
is a basis of $M^G$ as an $R$-module, there exist $\{b_i, 1\leq
i\leq n\}$ in $R$ such that
$$
dm=\sum_{1\leq i\leq n}da_i e_i=\sum_{1\leq i\leq n} b_i\sigma  e_i\ .
$$
It follows that for $1\leq i\leq n$ we have $da_i=b_i\sigma$  and thus
$b_i=dc_i$ with $c_i\in R$. Since $M$ is torsion
free we conclude that $ \sum_{1\leq i\leq n} a_ie_i=\sum_{1\leq
i\leq n} c_i\sigma e_i$ and hence that $m$ belongs to $M^G$. Since
$\ker\varphi''$ is torsion free and $\mbox{\rm coker}\varphi'$ 
is torsion, it follows
that $\ker\varphi''=\{ 0 \} = \mbox{\rm coker} \varphi'$.  
We then have therefore shown
that $\varphi'$ is an isomorphism.

 Let $\theta$ be the isomorphism
$
\theta: Hom_R(M_G,R) \rightarrow Hom_R(M^G,R)
$ 
induced by the isomorphism from $M^G$ onto $M_G$ described in
Lemma~\ref{fixedtrace}.
We want to describe $\theta\circ \varphi'$.
Let $x=\sigma m$ and $y=\sigma n $ be elements of $M^G$, then we
have the equalities
$
\theta\circ\varphi'(x)(y)=\varphi'(x)(n)=t(x,n)=t(m,y)=t^G(x,y) 
$.
We conclude that  $\varphi_{t^G}=\theta\circ\varphi'$ and
therefore that  $\varphi_{t^G}$  is an isomorphim as required.
\bigskip

\noindent For the
sake of completeness we conclude  this preparatory section  by
proving the following  well-known result (see \cite{McL} Cor.~3.3, p.~145
 and p.~196,
where the proof is carried out for the ``other'' action):

\begin{lemma} \label{ML}
Let $M$ and $N$ be left  $R[G]$-modules. Assume that $M$ and $N$
are both projective $R$-modules and that  either $M$ or $N$ is
projective as an $R[G]$-module. Then $M\otimes_R N$, endowed with the
diagonal action of $G$,  is a projective $R[G]$-module.
\end{lemma}

\Proof  Assume  $M$ to be $R[G]$-projective.  We must show  that
the functor $P \rightarrow Hom_{R[G]}(M\otimes_R N,P)$ from the
category of $R[G]$-modules into the category of $R$-modules is
exact. To this end  consider an arbitrary exact sequence of
$R[G]$- modules
$
0 \rightarrow  L'  \rightarrow  L  \rightarrow L'' \rightarrow  0
$.
Using the fact that $N$ is  projective as an $R$-module  we obtain
an exact sequence of $R$-modules
$$
0 \rightarrow  Hom_R(N,L')
 \rightarrow  Hom_R(N,L)  \rightarrow  Hom_R(N,L'') \rightarrow 0\ ,
$$
which is in fact readily seen  to be  an exact sequence  of
$R[G]$-modules. Using the fact  that $M$ is $R[G]$-projective, 
from the previous sequence we obtain a further  exact sequence of
$R$-modules
$$\begin{array}{rll}
0 \rightarrow  Hom_{R[G]}(M, Hom_R(N,L')) &
 \rightarrow  Hom_{R[G]}(M, Hom_R(N,L)) & \rightarrow\\
 & & \\
& \rightarrow  Hom_{R[G]}(M, Hom_R(N,L'')) & \rightarrow 0\ .
\end{array}$$

\noindent By the property of adjoint associativity of $Hom$ and $\otimes$
we obtain  for any
$R$-module $P$ a natural isomorphism of $R$-modules
$$ \psi_P:  Hom_R(M\otimes_R N,P) \simeq Hom_R(M,Hom_R(N,P))$$
$$f \mapsto (\psi_P(f):a \mapsto (b \mapsto f(a\otimes b)))\ .$$
Once again we check that $\psi_P$ is in fact an isomorphism of
$G$-modules and therefore induces an isomorphism of $R$-modules
$$\psi_P:  Hom_{R[G]}(M\otimes_R N,P) \simeq Hom_{R[G]}(M,Hom_R(N,P))\ .$$
Hence to conclude it suffices to replace the terms of the last exact sequence
by $Hom_{R[G]}(M\otimes_R N,P)$ for $P=L,L',L''$.
\bigskip

\subsection{Definition of the twist: the form}

Let $U$ be any open subscheme of $Y$, then $ R=O_Y(U)$ is an integral
 domain. Since the ramification of the cover is tame and $\pi$ is flat we
 know that $\pi _*({\cal D}^{-1/2}_{X/Y})(U)$ is a projective  $O_Y(U)[G]$-module, (see for instance section $4$). It  follows from Lemma~\ref{fixedtrace} and
 Lemma~\ref{ML} that $E_{\rho,X}(U)$
 is projective and therefore locally free over $O_Y(U)$.
 We then have proved that $E_{\rho,X}$ is a vector bundle over $O_Y$.
 We now endow the vector bundle $E_{\rho,X}$ with the form
 $$
 q_X={(q\otimes Tr_{X/Y})}^G \ .
 $$
\bigskip

\subsection{Proof of Theorem~\ref{itsabundle} (i) and examples of twists.} 

\noindent We  begin   this subsection by deducing Thm.~\ref{itsabundle}
 $(i)$ from the previous results.
For any affine open subscheme  $U$ of $Y$,
 we consider  the $O_Y(U)$-symmetric bundle $(M,t)$ with
 $M=(E\otimes  _{O_Y}\pi _*({\cal D}^{-1/2}_{X/Y}))(U)$
 and $t=(q\otimes Tr_{X/Y})$.
 We deduce from Prop.~\ref{CohAl} that $t^G$ defines a non-degenerate
 symmetric form on $E_{\rho,X}(U)$ and this   completes the proof of the
  first part  of theorem. The second part of the theorem will be 
proved in Sect.~4,  Lemma 4.3.  
 \bigskip

\noindent Next we give  some examples of twists of symmetric bundles.
\bigskip

\Examplenumb  Let $\mu$ denote the product form $\mu(x,y)=xy$ on 
$\pi _*({\cal D}^{-1/2}_{X/Y})$.
This  is a form in a generalised sense, in that it takes its values 
in the ideal $\pi _*({\cal D}^{-1}_{X/Y})$ instead of the ring $O_Y$. 
We  want to show that $q_{\rho ,X}$ is
nothing else than   the restriction to 
$E_{\rho ,X}={(E\otimes  _{O_Y}\pi _*({\cal D}^{-1/2}_{X/Y}))}^G$  
of the  form $t_X := q\otimes \mu $ on  
$(E\otimes  _{O_Y}\pi _*({\cal D}^{-1/2}_{X/Y}))$.   
\begin{prop}
$$
t_X(x,y) = q_{\rho , X}(x,y) \ .
$$
\end{prop}

\Proof In order to prove the equality it
suffices to prove  that $t_X$ and $q_{\rho ,X}$ coincide 
on $E_{\rho,X}(U)$ for any affine open $U$ in $Y$. 
Moreover any element of $E_{\rho, X}(U)$ can be obtained
as a sum of elements of the form
$\sigma (m\otimes n)$. So let $x=\sigma (m\otimes n)$ and
$y=\sigma (m'\otimes n')$. We have the equalities
$$
\begin{array}{rcl}
t_X(x,y) & = &
q_{\rho ,X}(\sum_{g \in G}(\rho(g) m\otimes gn),\sum_{h \in G}
\rho(h) m'\otimes hn')) \\
& & \\
& = & \sum_{g,h\in G}q(\rho(g)m,\rho(h)m')(gn)(hn')\ .
\end{array}
$$
Using the invariance of $q$ by $G$ this can be written
$$
t_X(x,y)=\sum_{g\in G}q(m,\rho(g)m'Tr_{X/Y}(n(gn'))) \ .
$$
We  now observe that  the right hand side of this equality can be
expressed as
$$
(q\otimes Tr_{X/Y})(m\otimes n,\sigma (m'\otimes n')) \ .
$$
Finally it follows from the very definition of $q_{\rho ,X}$ that this last
quantity is equal to $q_{\rho ,X}(x,y)$. Hence we have proved that
$t_X(x,y)=q_{\rho ,X}(x,y)$. 
\bigskip

\Examplenumb \label{subquotients} {\bf Subquotients.}
Here we show that the square root of the inverse different appears
as the twist of the standard
form on a permutation module by the natural orthogonal
representation attached to a Galois covering. This is to be compared with
the interpretation of the trace form of an \'etale algebra as a twist
of the standard form, which lies at the foundation of the original results
by Serre and Fr\"ohlich.
We keep the hypotheses and the notation of the general set-up as 
introduced in this section in $2$.a. We fix  a subgroup $H$ of $G$ and we let
$$ \lambda: X \rightarrow V:=X/H$$ 
denote the quotient map and 
$$\gamma: V \rightarrow Y$$ be the induced map, so that 
$\pi= \gamma \circ \lambda$. We note that $V$ being the  quotient of a 
normal scheme by a finite group is normal but not necessarily regular. 
Let $a$  run over  a left transversal of $H$ in $G$ and
let $E$ be the free bundle $O_Y[ G/H]$ which 
has for basis the left cosets $aH$ of $H$ in $G$.
We consider the symmetric $O_Y$-bundle $(E,q)$ 
where $q$ is the quadratic form on $E$ which has  $\{aH\}$ 
as an orthonormal basis. From now on we denote by $\bar a$ the coset $aH$. 
The left action of $G$ by permuting the cosets $\{\bar a\}$ 
extends to an orthogonal representation 
$$
\rho: G \rightarrow \Gamma (Y, {\bf O}(q)) \ .
$$
Our goal is to describe the twist of $(E,q)$ 
by the permutation representation $\rho$. 
  
\begin{prop} \label{Serre}
 \begin{itemize}
\item[i)] The twist of $(E,q)$ by $\rho$ is the symmetric bundle $(\pi _*({\cal D}^{-1/2}_{X/Y})^H,Tr_{V/Y})$ on $Y$. 

\item[ii)] Suppose that $\lambda$ is flat, then $(\pi _*({\cal D}^{-1/2}_{X/Y})^H,Tr_{V/Y})$ is the symmetric bundle $(\gamma_*({\cal D}^{-1/2}_{V/Y}), Tr_{V/Y})$.
\end{itemize}
\end{prop}
\Proof
 (i) We can reduce to an affine situation. 
For any point $y$ of $Y$ we denote by $U$ a sufficiently small  
open affine neighbourhood of $y$. Since the quadratic form $q$ is non-degenerate,  the adjoint map induces an isomorphism of $O_Y(U)G$-modules
$
E(U) \simeq Hom_{O_Y(U)}(E(U),O_Y(U))
$.
Tensoring by  ${\cal D}^{-1/2}_{X/Y}({\pi}^{-1}(U))$ and using [McL], V, Proposition 4.2, we obtain a further  isomorphism
$$
E(U) \otimes_{O_Y(U)} {\cal D}^{-1/2}_{X/Y}({\pi}^{-1}(U))\simeq 
Hom_ {O_Y(U)}(E(U), {\cal D}^{-1/2}_{X/Y}({\pi}^{-1}(U)))\ 
$$ 
that we denote by $\alpha$.  
We let  $G$ act diagonally on 
$E(U) \otimes_{O_Y(U)} {\cal D}^{-1/2}_{X/Y}({\pi}^{-1}(U))$ and by 
conjugation   on  $Hom_ {O_Y(U)}(E(U), {\cal D}^{-1/2}_{X/Y}({\pi}^{-1}(U)))$: 
$$
(g\cdot f)(u)=g(f(\rho(g^{-1})u), \ u \in E(U),\  g \in G \ .
$$ 
We observe that $\alpha$  respects the action of $G$. Therefore, by taking the fixed points of both sides, it induces by restriction an isomorphism:  
$$E_{\rho,X}(U) \simeq  Hom_ {O_Y(U)[G]}(E(U), {\cal D}^{-1/2}_{X/Y}({\pi}^{-1}(U)))\ .$$
We now consider the evaluation map 
$$
\begin{array}{rcl}
Hom_ {O_Y(U)[G]}(E(U), {\cal D}^{-1/2}_{X/Y}({\pi}^{-1}(U)))& \rightarrow &
{\cal D}^{-1/2}_{X/Y}({\pi}^{-1}(U))\\
 & & \\
f&\mapsto & f(\bar 1)\ . 
\end{array}
$$
For any $h \in H$ one has the equalities
$
h(f(\bar 1))=f(\bar h)=f(\bar 1)
$.
Hence we deduce that this  map takes its values in 
$ ({\cal D}^{-1/2}_{X/Y}({\pi}^{-1}(U)))^H$. Moreover we check
 easily that this is   an isomorphism. In summary we have constructed 
an isomorphism of 
$O_Y(U)$-modules 
$\theta: E_{\rho,X}(U) \rightarrow ({\cal D}^{-1/2}_{X/Y}({\pi}^{-1}(U)))^H$. 
Using Lemma 2.1, we  can now describe this isomorphism as follows.  
Let $x=\sigma (\bar a\otimes m) $ be an element of $E_{\rho,X}(U)$. 
The image of $x$ by $\alpha$  is the module homomorphism given by  
$(\bar b \rightarrow \sum_{g\in G}q(\bar {ga},\bar b)gm) \ .$ Therefore, 
we obtain that 
$$
\theta (x)=\alpha (x) (\bar 1)=\sum_{g\in G}q(\bar {ga},\bar 1)gm \ .
$$ 
It suffices to use the fact that $\{\bar {a}, a\in  G/H \}$ 
is an orthonormal basis to deduce  that 
$$\theta (x)=\theta (\sigma(\bar a\otimes m))=\sum_{h\in H}({ha^{-1}})m=Tr_{X/V}({a^{-1}}m) \ .$$
We now want to show that $\theta$ transports the 
form $q_{\rho,X}$ to the form $Tr_{V/Y}$. 
Let $x=\sigma(\bar a\otimes m) $ and  $y=\sigma(\bar b\otimes n) $ 
be elements of $E_{\rho,X}(U)$. From the very definition of $q_{\rho, X}$ 
we obtain the equality
$$
q_{\rho, X}(x,y)=\sum_{g\in G}q(\bar a,\bar {gb})Tr_{X/Y}(m(gn))\ .
$$
Using once again the fact that $\{\bar {a}, a\in  G/H \}$ is 
an orthonormal basis, we deduce from the above that 
$$
q_{\rho, X}(x,y)=\sum_{h\in H}Tr_{X/Y}(m((ah{b^{-1}})n))=\sum_{h\in H}Tr_{X/Y}(((h^{-1}{a^{-1}})m)({b^{-1}}n))\ . 
$$
Therefore we have proved that 
$$
q_{\rho, X}(x,y)=Tr_{V/Y}(\sum_{h,k\in H}h({a^{-1}}m)k({b^{-1}}n))=Tr_{V/Y}(\theta (x)\theta (y)) \ .
$$
We conclude as required that $\theta$ is an isomorphism of symmetric $O_Y(U)$-bundles
$$ 
(E_{\rho,X}(U),q_{\rho,X}(U))\simeq  
(({\cal D}^{-1/2}_{X/Y}({\pi}^{-1}(U)))^H, Tr_{V/Y}) \ .
$$
(ii) In order to complete the proof of the proposition 
it suffices to show that, if $\lambda$ is flat, then:
$$({\cal D}^{-1/2}_{X/Y}({\pi}^{-1}(U)))^H={\cal D}^{-1/2}_{V/Y}({\gamma}^{-1}(U)) \ . $$
>From the transitivity formula for the differents this last equality
reduces us to showing  that 
$$({\cal D}^{-1/2}_{X/V}({\pi}^{-1}(U))^H=O_V(\gamma^{-1}(U)) \ .$$
Since  the cover $X \rightarrow V$  is tame and $\lambda$ is flat,  it follows that  ${\cal D}^{-1/2}_{X/V}({\pi}^{-1}(U))$ is a projective  $O_V(\gamma^{-1}(U))[H]$-module. Once again from Lemma 2.1 we deduce that
$$ ({\cal D}^{-1/2}_{X/V}({\pi}^{-1}(U)))^H=Tr_{X/V}( {\cal D}^{-1/2}_{X/V}({\pi}^{-1}(U)))\ .$$
 From \cite{SerreCL}, III, Prop.~7, we  obtain that 
$$Tr_{X/V}( {\cal D}^{-1/2}_{X/V}({\pi}^{-1}(U)))=O_V(\gamma^{-1}(U)) \ .$$
The required equality now follows.
\bigskip

\Examplenumb 
Let $\pi_1(Y)^{t,o}$ be the tame odd fundamental group of $Y$. 
Any continuous homomorphism $\rho: \pi_1(Y)^{t,o}  \rightarrow  
{\bf Z}/2{\bf Z }$ can be seen as an orthogonal representation of a 
finite quotient of $\pi_1(Y)^{t,o}$ into $\Gamma (Y,{\bf O}(q))$ 
where $(E,q)$ is the standard rank one  symmetric bundle $(O_Y,x^2)$. 
Therefore, let $X/Y$ be the \'etale cover defined by $\rho$, then the 
map  $\rho \mapsto (E_{X,\rho},q_{X,\rho})$ induces a  canonical 
map
$$
can: H^1(\pi_1(Y)^{t,o}, {\bf Z}/2{\bf Z}) \rightarrow 
H^1(Y,{\bf Z}/2{\bf Z }) \ .
$$  
 
\bigskip

\section{The \'etale case}

\noindent In this section we study  the case  
where $X$ is a torsor of the constant
group scheme $G$ over $Y$, so  that $X$ is an \'etale covering
of $Y$. We consider, as earlier, a symmetric
bundle $(E,q)$, an orthogonal representation $\rho:G\rightarrow
\Gamma (Y,{\bf O}(q))$ and the twist $(E_{\rho ,X},q_{\rho ,X})$ of $(E,q)$ by
$\rho$  which,  as observed previously,  coincides with  the
$O_Y$-vector bundle $(\pi_*\pi^*E)^G$.

 Our aim in this section is to prove, as announced in  Thm.~0.3,
 comparison formul{\ae} between the first and the second
 Hasse-Witt invariants of $(E,q)$ and $(E_{\rho ,X},q_{\rho ,X})$ 
which generalise
 those obtained by Fr\"ohlich, [F], Thm.~2 and Thm.~3, in
 the case of fields extensions.

 On the one hand, following Milne \cite{Mi}, III, Proposition 4.6,
 we note that the isomorphism class of  $X$, considered as a sheaf
 torsor for $G$,  defines an element $c(X)$ in the cohomology set
 $H^1(Y,G)$. On the other hand we show that $(E_{\rho,X},q_{\rho,X})$ is a
 twisted form of $(E,q)$; therefore its isometry class defines an element
 denoted by  $[E_{\rho ,X},q_{\rho ,X}]$ in  $ {H}^1(Y,{\bf O}(q))$. 
Finally the morphism
 $\rho$ induces a natural map $\rho_*: { H}^1(Y,G)
 \rightarrow { H}^1(Y,{\bf O}(q))$. As announced in the introduction we 
establish the following connection between these objects.

\begin{prop} \label{torsor}  In  $ { H}^1(Y,{\bf O}(q))$ one
has the equality
 $$
\rho_*(c(X))=[E_{\rho ,X},q_{\rho ,X}]\ .
$$
\end{prop}

\Proof For any sheaf of groups  in  the \'etale topology, $\cal F$
on $Y$,   the set $ { H}^1(Y,\cal F)$ is defined to be
$\displaystyle  \lim_{\rightarrow} {H}^1(\cal U,\cal F)$
where the limit is taken over all \'etale coverings $\cal U$ of $Y$.
The strategy of the proof is to show that, for any ``sufficiently fine''
\'etale covering  $\cal U$ $ =(U_i \rightarrow Y)_{i\in I}$
 which trivializes $X$ as a torsor for $G$,
 we obtain an isometry of symmetric bundles from $(E_{\rho,X}\times_Y U_i, q_{\rho,X})$
 onto $(E\times_Y U_i,q)$. For such coverings $\cal U$ we will
 first obtain a $1$-cocycle representing $c(X)$ with values in $G$, then
 from this we obtain  a $1$-cocycle representing $[E_{\rho ,X},q_{\rho ,X}]$ 
which takes values in ${\bf O}(q)$.

Let ${\cal U} = (U_i)_{i\in I}$ be a sufficiently fine, affine, \'etale
cover which trivialises $X$ as a $G$-torsor. We denote respectively  by $O_Y(U_i)$
and $O_X(U_i)$ the  global sections of $U_i$ and $X\times_Y U_i$.
Now $X\times_Y U_i \rightarrow U_i$  is finite and $U_i$ affine,
thus $X\times_Y U_i$ is affine. Therefore for any $i$ the
isomorphism of $U_i$-schemes $G_{(U_i)} \simeq X\times_Y U_i$  is
induced from an $O_Y(U_i)$-$G$ isomorphism of algebras
$$
\Phi_i: O_X(U_i) \rightarrow Map(G, O_Y(U_i)) \ .
$$
Furthermore we note that $g_{ij}=\Phi_j{\Phi_i}^{-1}$ is the
$1$-cocycle representing $c(X)$, (\cite{Mi}, III, section $4$).

We again  denote by $Tr_{X/Y}$ the form induced by the trace on
$O_X(U_i)$. For any elements $x$ and $y$  of $O_X(U_i)$, the trace
$Tr_{X/Y}(xy)$ belongs to $O_Y(U_i)$. Therefore, since $\Phi_i$ is
an isomorphism of $O_Y(U_i)$-algebras, we will have
$$  
Tr_{X/Y}(xy)=\Phi_i(Tr_{X/Y}(xy))=\Phi_i(Tr_{X/Y}(xy))(1) \ .
$$
Using now that $\Phi_i$ is $G$-equivariant we  deduce from the
previous equalities that
$$  Tr_{X/Y}(xy)=\sum_{g\in G} \Phi_i(x)(g))\Phi_i(y)(g) \ .$$
Denoting by $\mu_G$ the standard $G$-invariant form on $Map(G, O_Y(U_i))$ given
by
$$\mu_G(f,f')=\sum_{g\in G} f(g)f'(g) \ ,$$
we have then proved that $\Phi_i$ is an equivariant  isometry
between the $O_Y(U_i)-G$ symmetric bundles $(O_X(U_i), Tr_{X/Y})$
and $( Map(G, O_Y(U_i)),\mu_G)$.

Let $(E(U_i),q)$ be the $O_Y(U_i)$ symmetric bundle defined by
considering the global sections of the inverse image of $(E,q)$ by
the morphism $U_i \rightarrow Y$. After tensoring over $O_Y(U_i)$
and taking fixed points by $G$ we deduce from $\Phi_i$ an isometry
that we again  denote by $\Phi_i$
$$ \Phi_i:( E_{\rho,X}(U_i),q_{\rho,X}) \rightarrow (E(U_i)\otimes_{O_Y(U_i)} Map(G, O_Y(U_i)))^G, t^G) \ ,$$
where $t$ is the form $q\otimes \mu_G$ and $t^G$ is obtained from
$t$ by following the recipe described in Sect.~2. We now let
$\nu_i$ be the morphism of $O_Y(U_i)$-modules induced by
$$ 
\begin{array}{ccl}
\nu_i:E(U_i)\otimes_{O_Y(U_i)} Map(G, O_Y(U_i)) & \longrightarrow & 
Map(G, E(U_i))\\
 & & \\
e\otimes f & \mapsto & (g \mapsto f(g)(\rho(g)e)) \ .
\end{array}
$$
It is easy to check that $\nu_i$ is an isomorphism (use for
instance the fact that \newline
$Map(G, O_Y(U_i))$ is a free  $O_Y(U_i))[G]$-module of rank $1$
with basis $l$, where $l$  is defined by $l(g)=1$ if
$g=1$ and $0$ otherwise). The group $G$ acts diagonally on the
left hand side while on the right hand side it acts by
$uf:g \mapsto  f(gu)$.
It follows from the definitions of these actions that $\nu_i$ is a
$G$-isomorphism and thus induces an isomorphism, again  denoted
$\nu_i$
$$ \nu_i:(E(U_i)\otimes_{O_Y(U_i)} Map(G, O_Y(U_i)))^G \rightarrow  
Map(G, E(U_i))^G \ .$$ We now observe that  the map $f \mapsto f(1)$ 
is clearly an isomorphism from $ Map(G, E(U_i))^G$ onto $E(U_i)$. Hence, 
finally composing this map with $\nu_i$ we have defined an isomorphism of 
$O_Y(U_i)$-modules
$$
\gamma_i : (E(U_i)\otimes_{O_Y(U_i)} Map(G, O_Y(U_i)))^G 
\rightarrow E(U_i) \ .
$$
We now want to describe in some detail the map $\gamma_i$. We note that,  
since the set $\{gl,\ g\in G\}$ is a free basis of the $O_Y(U_i)$-module 
$Map(G, O_Y(U_i))$, every element  of  $ (E(U_i)\otimes_{O_Y(U_i)} 
Map(G, O_Y(U_i))$ can be written as a sum $\sum_{g\in G} a_g\otimes gl$ and 
therefore every element of the subgroup of the fixed points by $G$ can 
be written as a sum $\sum_{g\in G} \rho(g)e\otimes gl$. 
Let  $x$ be an element of    
$ (E(U_i)\otimes_{O_Y(U_i)} 
Map(G, O_Y(U_i))^G$ written as such a sum, we then have:
$$ 
\gamma_i(x)=\nu_i(x)(1)=\sum_{g\in G}\nu_i( \rho(g)e\otimes gl)(1)=
 \sum_{g\in G} (gl)(1)\rho(g)e=e \ .
$$
We now also consider  
$y=\sum_{g\in G} \rho(g)e'\otimes gl$. 
Since $x=\sigma (e\otimes l)$ and $y=\sigma (e'\otimes l)$  
it follows from the very definition of $t^G$, prior to 
Prop.~\ref{CohAl},  that
$$t^G(x,y)=t(x,e'\otimes l)=\sum_{g\in G}q(\rho(g)e,e')\mu_G(gl,l) \ .$$
Since $\{gl, g\in G\}$ is an orthonormal basis of $Map(G, O_Y)$,
the right hand side of the last equality is equal to $q(e,e')$. We
conclude that $\gamma_i$ is an isometry and therefore that
${\Phi_i}^{-1}{\gamma_i}^{-1}$ is an isometry $\theta_i$
$$\theta_i:(E(U_i),q) \rightarrow (E_{\rho,X}(U_i),q_{\rho,X})\ .$$
We now must evaluate ${\theta_i}^{-1}\theta_j$ in order to obtain
a $1$-cocycle representing $(E_{\rho,X},q_{\rho,X})$ as a twisted form of
$(E,q)$. Starting with $e \in E(U_{ij})$ we obtain that
$$
(\gamma_i\Phi_i{\Phi_j}^{-1}{\gamma_j}^{-1})(e)=\gamma_i(\sum_{g\in
G}\rho(g)e\otimes (g{g_{ij}^{-1}})l) = \gamma_i(\sum_{u\in
G}(\rho(ug_{ij})e)\otimes ul) = \rho (g_{ij})e \ .
$$
This concludes the proof of the proposition.
\bigskip

\noindent {\bf Remark.}
The twist of a  symmetric bundle $(E,q)$  by an  orthogonal representation $\rho$ is always a twisted form of $(E,q)$ and therefore is given by a class in $H^1(Y, {\bf O}(q))$. In the \'etale situation  Prop.~\ref{torsor} tells  us precisely that this class is the image by $\rho_*$ of the class defining $X$ as
a torsor for $G$. Therefore the determination of the Hasse-Witt invariants of the twisted bundle will be obtained as an application of Theorem 0.2.   

\smallskip

\subsection{Proof of Theorem~\ref{etale}}

\noindent We now return to the situation considered in Prop.~\ref{torsor} and
we generalise  Theorems 2 and 3  in [F], from field extensions to \'etale 
covers. We start by recalling some  notation and definitions. 
Let $\pi_1(Y)$ be the
fundamental group of $Y$ based at some chosen  geometric point.
We consider a representation $\rho :\pi_1(Y) \rightarrow \Gamma
(Y,{\bf O}(q))$, where $(E,q)$ is a symmetric $Y$-bundle. We assume
$\rho$ to have an open kernel $N$. Then $N$ defines a finite
Galois \'etale cover $X/Y$ with Galois group $G=\pi_1(Y)/N$. The
cohomology class $c(X)$ of $H^1(Y,G)$ defined  by $X$, considered
as a $G$-torsor,  only depends on $\rho$ and therefore will be
denoted by $c(\rho)$. The representation $\rho$ factorises into a
homomorphism $\rho :G \rightarrow \Gamma(Y,{\bf O}(q))$ which in turn
induces $\rho_* : H^1(Y,G) \rightarrow H^1(Y,{\bf O}(q))$
  and $\rho_{*\bar K}: H^1(Y,G) \rightarrow H^1(Y,{\bf O}(q)(\bar K))$
  with $\bar K$ denoting the separable closure of the residue field at
  some point. 
We also consider the two exact sequences 
$$ 
1\rightarrow {\bf Z}/2{\bf Z} \rightarrow {\widetilde {\bf O}}(q)
\rightarrow {\bf O}(q) \rightarrow 1 \ ,
$$ 
and
$$
 1\rightarrow {\bf Z}/2{\bf Z} \rightarrow {\widetilde {\bf O}}(q)(\bar K)
\rightarrow {\bf O}(q)(\bar K) \rightarrow 1 \ ,
$$
where the first  is an exact sequence of \'etale sheaves of groups and the
 second is the exact sequence of the $\bar K$-points of the first. These
sequences induce the boundary maps
$$
\delta^2: H^1(Y, O(q)) \rightarrow H^2(Y, {\bf Z}/2{\bf Z})
$$
and
$$
\delta^2_{\bar K}: H^1(Y, O(q)(\overline K))
\rightarrow H^2(Y, {\bf Z}/2{\bf Z}) \ .
$$

\Definition \label{sw2-sp2}  
The {\it first Stiefel-Whitney class} of $\rho$ is defined
to be  
$$w_1(\rho)=\delta^1(\rho_*(c(\rho)) \ .
$$
The {\it second  Stiefel-Whitney class} of
$\rho $ is defined by 
$$
w_2(\rho)= (\delta^2_{\bar K}\rho_{*\bar K})(c(\rho)) \ .
$$ 
The {\it spinor class} is defined to be the difference
$$
sp_2(\rho)=(\delta^2\rho_*)(c(\rho))-(\delta^2_{\bar K}\rho_{*\bar K})(c(\rho))
\ .
$$

\bigskip

\noindent {\bf Remark.} The notion of the spinor class was first
introduced by Fr\"ohlich for the case of field extensions.
Fr\"ohlich's initial definition was generalised by Khan in a
geometric context. In fact it is not immediately clear that Kahn's
definition coincides with  ours. Nevertheless one has to observe
that in the case of field extensions, when $Y=\mbox{Spec} (K)$ and $K$ is a
field of characteristic different from $2$, then  T. Saito has proved in 
[Sa1] Lemma 3, that the spinor class we consider here is indeed equal to
Fr\"ohlich's original one.

\smallskip

The proof of Thm.~0.4 is now an easy consequence of
Prop.~\ref{torsor} and Thm.~\ref{general}. Let us denote by $\alpha$ the
cohomology class $\rho_*(c(\rho))$. It follows from Prop.~\ref{torsor}
that $(E_{\rho , X},q_{\rho ,X})$ can be taken as a representative 
of the isometry
class of twisted forms of $(E,q)$ whose image is $\alpha$.
Therefore it follows from  Thm.~\ref{general}  that
$w_1(E_{\rho , X})=w_1(E)+\delta^1(\alpha)$ and
$w_2(E_{\rho , X})=w_2(E)+w_1(E)\delta^1(\alpha) +\delta^2(\alpha)$. From
the above  definitions  we deduce that $
\delta^1(\alpha)=w_1(\rho)$ and that
$\delta^2(\alpha)=w_2(\rho)+sp_2(\rho)$. It now suffices to use these equalities  in order to
obtain Thm.~\ref{etale}.
\bigskip

\section{The tame case}

\noindent In this section we return to the general tame situation as 
described in Sect.~2.a.  
We let $X$ be  a connected projective regular scheme
equipped with a tame action of the  finite group $G$ and let $Y$ denote
$X/G$. We assume that $\pi : X \rightarrow Y$ is flat and thus that $Y$ 
is regular. We consider a symmetric $Y$-bundle $(E,q)$,
an orthogonal representation $\rho :G \rightarrow \Gamma (Y, {\bf O}(q))$ 
and we let
$(E_{\rho ,X},q_{\rho ,X})$ be the twist of $(E,q)$ by $\rho$ as defined
previously.  As in the \'etale case, our aim is to compare the
Hasse-Witt invariants of $(E,q)$ and $(E_{\rho ,X},q_{\rho ,X})$. 
As previously, since
$(E_{\rho ,X},q_{\rho ,X})$ is a twist 
of $(E,q)$, it defines an element $\alpha$
in the set of  cohomology classes  $H^1(Y, {\bf O}(q))$ and  we
can therefore apply the results of Thm.~0.2.  Nevertheless we cannot yet
give a description of the class $\alpha$  as explicit as in
the \'etale case.  In order to obtain an explicit comparison
formula,  we shall therefore associate to our given tame cover an
\'etale cover to which  we can directly  apply our previous results. 
As explained in the introduction, 
the difference between the formula for the
tame cover  we started with and that for the \'etale cover that 
we construct will
be reflected in the  appearance  of a new class
depending on the decomposition of the representation $\rho$ when
restricted to the inertia groups of the generic points of the
irreducible components of the branch locus of the covering $X/Y$. Two 
technical notions will play an important role here: that of a metabolic bundle
and that of normalisation along a divisor.
\bigskip

\subsection{Metabolic bundles and their twists}

Hyperbolic forms are
metabolic and the greater generality of the latter notion is the
correct one when dealing with forms over  non-affine schemes. 
We recall the definition.

Let $E$ be a vector bundle over $Y$. A sub-$O_Y$-module $V$ of $E$
is a {\it sub-bundle} of $E$ if it is {\it locally} a direct
summand, {\it i.e.} for any $y$ in $Y$, there is an open $Z$
containing $y$ such that $V|_Z$ has a direct summand in $E|_Z$. If
$V$ is a sub-bundle of $E$, then $V$ and the quotient $E/V$ are
both vector bundles. Let now $(E,q)$ be a vector bundle endowed
with a quadratic form. For a  sub-module $i : V\subset E$ one
defines an orthogonal complement, which is the sub-$O_Y$-module
$V^{\bot}$, whose sections over the open $Z$ consist of those
sections of $E$ which are orthogonal to all sections of $V$ over
any open subset of $Z$. Alternatively:
$$
V^{\bot} = \ker (E\stackrel{\varphi_q}{\rightarrow}E^{\vee}
\stackrel{i^{\vee}}{\rightarrow} V^{\vee}) \ .
$$
If furthermore $V$ is a sub-bundle of $E$, then $i^{\vee}$ is an
epimorphism. Assume now that  $(E,q)$ is a symmetric bundle, then
$\varphi _q$ is by definition an isomorphism and therefore we have
an isomorphism
$$
\alpha : E/V^{\bot}\cong V^{\vee} \ ,
$$
$E/V^{\bot}$ is locally free and  $V^{\bot}$ is also a sub-bundle.
There also is an isomorphism
$$
\beta : V^{\bot}\cong (E/V)^{\vee} \ .
$$
(The module  $(E/V)^{\vee}$ can be identified with the
sub-$O_Y$-module of $E$ whose sections over $Z\subset Y$ are the
linear forms $\lambda : E|_Z\rightarrow O_Z$ which vanish on $V$,
so $V^{\bot}=\varphi _{q} ^{-1}((E/V)^{\vee})$.)

A sub-bundle $V$ of a symmetric bundle  $(E,q)$ is a {\it totally
isotropic sub-bundle} (also called a {\it sub-lagrangian}) if
$V\subset V^{\bot}$. If $V$ is a sub-lagrangian of $(E,q)$, then
$V^{\bot}/V$ is a sub-bundle of $E/V$ and the form on $V^{\bot}/V$
obtained by reducing $q$ modulo $V$ is non-degenerate. A
sub-bundle $V$ of   $(E,q)$ is called a {\it lagrangian} if it is
such that $V=V^{\bot}$. The symmetric bundle  $(E,q)$ is called
{\it metabolic} if it contains a lagrangian. If $V$ is a
lagrangian in $(E,q)$, then $\mbox{rank}(E)=2\cdot\mbox{rank}(V)$
and $V$ is in a sense a maximal totally isotropic sub-bundle. One
can observe  that $V$ is a lagrangian in $(E,q)$ if and only if
one has a commutative diagram
$$
\begin{array}{ccccccccc}
0& \rightarrow & V & \rightarrow & E & \rightarrow &V^{\vee} &\rightarrow & 0\\
& &id  \downarrow \ \ \ \ & & \downarrow  \varphi _q & &\ \ \
\downarrow id
& &  \\
0& \rightarrow & V\cong V^{\vee\vee}  & \rightarrow & E^{\vee} &
\rightarrow & V^{\vee} &\rightarrow &0 \ .
\end{array} 
$$
That is to say: a metabolic form is given by a symmetric self-dual short
exact sequence (see [Kne],  Chapt.~3). This notion can be generalized to
symmetric complexes (see Appendix, \cite{Ba1} and \cite{Ba2}).
Metabolic bundles are trivial in the Witt group
$W\left(  Y\right)  $, but the converse does not in general hold 
(see \cite{O2} 
Ex.~2.10 for an example). 

\bigskip

\Example {\bf The twist of a metabolic bundle is metabolic.}
Let $(E,q)$ be the underlying symmetric $Y$-bundle of an 
orthogonal representation $\rho:G \rightarrow \Gamma (Y,{\bf O}(q))$, 
where $X$ is a scheme endowed with a tame action by the finite group 
$G$ and $Y=X/G$. We shall say that $\rho$ is  a metabolic representation when  
$(E,q)$ is metabolic with a $G$-invariant lagrangian $V$.  
Our goal is to prove 
that if $(E,q)$ and $\rho$ are metabolic then 
$(E_{\rho,X},q_{\rho,X})$ is also metabolic.

Let $V$ be a lagrangian of $(E,q)$. We then  have an exact sequence of 
$G$-modules
$$
0\rightarrow  V \rightarrow  E  \rightarrow V^{\vee} \rightarrow  0 \ .
$$
Since the cover $\pi:X \rightarrow Y$ is tame and $\pi$ is flat, 
$\pi_*({\cal D}_{X/Y}^{-1/2})$ is a locally projective $O_Y[G]$-module. 
Therefore it follows from Lemma 2.3 that the modules 
$V\otimes_{O_Y}\pi_*({\cal D}_{X/Y}^{-1/2})$,  
 $E\otimes_{O_Y}\pi_*({\cal D}_{X/Y}^{-1/2})$ and 
$V^{\vee}\otimes_{O_Y}\pi_*({\cal D}_{X/Y}^{-1/2})$ are all locally 
projective  $O_Y[G]$-modules. So,  tensoring and taking $G$-fixed 
points, affords  a new exact 
sequence of locally free $O_Y[G]$-modules
 $$
0 \rightarrow  (V\otimes_{O_Y}\pi_*({\cal D}_{X/Y}^{-1/2}))^G  \rightarrow  
E_{\rho,X}  \rightarrow (V^{\vee}\otimes_{O_Y}\pi_*({\cal D}_{X/Y}^{-1/2}))^G 
\rightarrow  0\ ,
$$  
(see the proof of Prop.~4.4 for details). 
It  is clear from  the definition of $q_{\rho,X}$ that its restriction  
to $(V\otimes_{O_Y}\pi_*({\cal D}_{X/Y}^{-1/2}))^G$ is the null form since 
$V$ is a lagrangian for $q$. In order to prove  that $E_{\rho,X}$ is 
metabolic with $(V\otimes_{O_Y}\pi_*({\cal D}_{X/Y}^{-1/2}))^G$ as a 
lagrangian it suffices to show that the rank of  $E_{\rho,X}$ is 
twice the rank of  $(V\otimes_{O_Y}\pi_*({\cal D}_{X/Y}^{-1/2}))^G$. 
This will follow at once from the $O_Y$-module isomorphism
$$ 
((V\otimes_{O_Y}\pi_*({\cal D}_{X/Y}^{-1/2}))^G)^{\vee} 
\simeq  (V^{\vee}\otimes_{O_Y}\pi_*({\cal D}_{X/Y}^{-1/2}))^G \ .
$$
This last isomorphism can be proved as follows. Since $\pi_*({\cal D}_{X/Y}^{-1/2})$ is unimodular it is self-dual and therefore
 $$
(V^{\vee}\otimes_{O_Y}\pi_*({\cal D}_{X/Y}^{-1/2})) 
\simeq  (V\otimes_{O_Y}\pi_*({\cal D}_{X/Y}^{-1/2}))^{\vee} \ .
$$
Hence we deduce that 
$
(V^{\vee}\otimes_{O_Y}\pi_*({\cal D}_{X/Y}^{-1/2}))^G
$
is isomorphic to
$$(V\otimes_{O_Y}\pi_*({\cal D}_{X/Y}^{-1/2})^{\vee})^G=
Hom_{O_Y} ((V\otimes_{O_Y}\pi_*({\cal D}_{X/Y}^{-1/2}))_G,O_Y) \ .
$$
Since  we know that   $(V\otimes_{O_Y}\pi_*({\cal D}_{X/Y}^{-1/2})$ 
is a locally projective $O_Y[G]$-module,
 it follows from \noindent Lemma 2.1 $(ii)$  that  $(V\otimes_{O_Y}\pi_*({\cal D}_{X/Y}^{-1/2}))_G$ is isomorphic to $ (V\otimes_{O_Y}\pi_*({\cal D}_{X/Y}^{-1/2}))^G$. We then conclude that 
$$Hom_{O_Y} ((V\otimes_{O_Y}\pi_*({\cal D}_{X/Y}^{-1/2}))_G,O_Y)  \simeq 
 ((V\otimes_{O_Y}\pi_*({\cal D}_{X/Y}^{-1/2})^G)^{\vee} \ , $$ as required.  
  
\subsection{The main lemma}
Here we state a result of \cite{EKV}, Prop.~5.5,  which is one of
the main tools in calculating Hasse-Witt invariants. It provides a
precise relationship  between the Hasse-Witt invariants of a metabolic
form and the Chern classes of a lagrangian. 
 In  \cite{CNET2} we have shown  the relevance of
this lemma when dealing with symmetric complexes (see the Appendix). 
We recall that we
have associated to any  bundle $V$ of finite rank  over $Y$ an element
$d_t(V)$ of $H^*(Y, {\bf Z}/2{\bf Z})$ defined prior to Thm.~0.8 in the 
introduction. 

 \begin{lemma}
Let $(E,q)$ be a metabolic form with lagrangian $V$. Then
$$ w_t(E) = d_t(V) \ .$$
In particular
$$1+w_1(E)t+w_2(E)t^2=1+n(-1)t+\left( c_1(V)+\left( \! \! \begin{array}{c} n \\ 2
                      \end{array} \! \! \right)(-1)\cup (-1)\right)t^2 \ .$$
\end{lemma}

Let us now assume $Y$ to be  irreducible  with generic point $\eta$.
In our set-up we shall  construct a number of   metabolic bundles from symmetric
bundles over $Y$  which are isometric when restricted to the
generic fiber of $Y$. The essential point here is that  if $(E,q_E)$ and $(F,q_F)$ are
defined over $Y$ and agree on the generic fiber
$$
E|_{\eta } = F|_{\eta } \ ,
$$
then, under suitable assumptions,  $(E\bot F, q_E\bot -q_F)$ is
metabolic. A natural approach is then to consider  the sub-sheaf
 ${\cal G}$ of $E|_{\eta} = F|_{\eta}$ defined as
$$
{\cal G} := < e-f | e\in E\ , f\in F > \ ,
$$
and the exact sequence
$$
0\rightarrow E\cap F \rightarrow E\bot F \rightarrow {\cal
G}\rightarrow 0 \ ,
$$
where the maps are obtained by restricting to $E\bot F$ the maps
defined at the level of the generic fibers given by the diagonal
map and the map sending $(x,y)$ to $(x-y)/2$. Then,  if $ {\cal
G}$ is locally free, it follows that $E\cap F$ is a sub-bundle of
$E\bot F$, which is a lagrangian and
$$
E\bot F / E\cap F \simeq (E\cap F)^{\vee} \simeq {\cal G} \ .
$$
\begin{lemma} \label{lemma-main}
Let $(E,q_E)$ and $(F,q_F)$ be non-degenerate forms which are
isometric when restricted to $\eta$.
\begin{itemize}
\item[a)] $w_1(E)=w_1(F)$ and $c_1(E)=c_1(F)$.
\item[b)]
Consider the exact sequence
$$
0\rightarrow E\cap F \rightarrow E\bot F \rightarrow {\cal
G}\rightarrow 0 \ .
$$
Assume that ${\cal G}$ (and therefore $E\cap F$) is locally free over $O_Y$, and  that
$(E\bot F, q_E\bot -q_F)$ is metabolic with lagrangian $E\cap
F={\cal G}^{\vee}$. Then
$$
w_t(E,q_E)\cdot w_t(F,-q_F)=d_t({\cal G})  \ .
$$
\item[c)] With the assumptions and the notations of (b), 
$$w_2(E,q_E)-w_2(F,q_F)=c_1({\cal G}) -c_1(E) \ .$$
 In particular
this sum belongs to the image of the Picard group modulo $2$ in
$H^2(Y, {\bf Z}/2{\bf Z})$. 
\end{itemize}
\end{lemma}
\Proof Part $(a)$ of the Lemma is proved in \cite{EKV}, Lemma 4.8; 
$(b)$ is an immediate  corollary of  Lemma 4.1; and   $(c)$ can also 
be deduced from Lemma 4.1, (see  \cite{EKV}, Cor.~6.3).  
\medskip

\noindent {\bf Remark.}
The first Chern class $c_1(L)$ of a vector bundle $L$ over $Y$ is
obtained as the image in  $\mbox{\rm Pic}(Y)$ of the class defined by the line
bundle $\mbox{\rm det}(L)$ under the boundary map associated to the Kummer
exact sequence deduced from the multiplication by $2$ map on 
${\bf G}_m$.

\subsection{Normalisation along a divisor}

\noindent The  process of normalisation along a branch divisor  
was introduced and studied in \cite{CNET1}. Let $G_2$ be  a $2$-Sylow subgroup of $G$. We write
$Z={ X}/G_2$ and we let  $T$ be  the normalisation of the fiber
product $T'=Z\times _YX$. So we have the  diagram
$$
\begin{array}{ccccl}
T & \longrightarrow & T'= Z\times _YX & \stackrel{\phi'}{\rightarrow} & X \\
 \ & & & & \\
 &\pi _Z \searrow & \downarrow & & \downarrow \pi \\
 \ & & & & \\
 &  & Z & \stackrel{\phi}{\rightarrow} & Y
\end{array}
$$
We have proved in \cite{CNET1}, Thm.~2.2 that $T$ is regular and that the map 
$\pi _Z$ is \'etale. In {\it loc. cit}, Sect.~3.3 we showed how to decompose
the normalisation map $T\rightarrow T'$ into a sequence
$$
T=T^{(m)}\rightarrow T^{(m-1)}\rightarrow \cdots \rightarrow
T^{(0)}=T' \ ,
$$
of flat $Z$-covers, with each $\pi ^{(h)} : T^{(h)}\rightarrow Z$
having the property that 
${\cal D}^{-1/2}_{T^{(h)}/Z}$ is a well defined $T^{(h)}$-vector bundle.
We set $\Lambda ^{(h)} = \pi _{*}^{(h)}({\cal D}^{-1/2}_{T^{(h)}/Z})$ 
and define $I^{(h)} = \Lambda ^{(h)} \cap \Lambda ^{(h+1)}$ and
${\cal G}^{(h)} = (\Lambda ^{(h)} +  \Lambda ^{(h+1)})/I^{(h)}$.
Then for $0\leq h\leq m-1$ we have short exact 
sequences of locally free
$O_Z$-modules
$$
0\rightarrow I^{(h)}\rightarrow \Lambda ^{(h)}\oplus \Lambda
^{(h+1)} \rightarrow {\cal G}^{(h)}\rightarrow 0 \ .
$$
The $\Lambda ^{(h)}$ all coincide on the generic fiber and we have proved
that $I^{(h)}$ and ${\cal G}^{(h)}$ are locally free $O_Z$-modules.
We then deduced that for $1\leq h \leq m-1$, 
$(\Lambda ^{(h)}, (-1)^hTr_{T^{(h)}/Z}) \perp 
( \Lambda ^{(h+1)}, (-1)^{(h+1)}Tr_{T^{(h+1)}/Z})$ is a metabolic bundle
which satisfies the hypotheses of Lemma 4.2.

\subsection{Proof of Theorems 0.6 (ii), 0.7  (i) and \ref{big-main-lemma}}

\noindent  
As explained in the introduction,   our strategy  consists of 
considering

\noindent $w_k(\phi ^{*}(E)_{\phi ^{*}(\rho),T'})-w_k(\phi ^{*}(E))$ as a sum of
two terms, namely:
$$
(w_k(\phi ^{*}(E)_{\phi ^{*}(\rho),T'})-w_k(\phi ^{*}(E)_{\phi ^{*}(\rho),T})) 
+
(w_k(\phi ^{*}(E)_{\phi ^{*}(\rho),T})-w_k(\phi ^{*}(E))
 \ .
$$
The cover $T/Z$ is \'etale and so it is a $G$-torsor. Hence
the second term is known by Thm.~\ref{etale} for
$k=1, 2$. Therefore the main part of this
section will be devoted to computing the first term.

Our first aim is to prove that our twisting process behaves well under
flat base change (see Thm.~0.6 (ii)).

 \begin{lemma} For any orthogonal representation 
$\rho :G \rightarrow \Gamma (Y, {\bf O}(E))$
 $$
 \phi ^{*}(E)_{\phi ^{*}(\rho), T'}=\phi ^{*}(E_{\rho ,X}) \ .
 $$
\end{lemma}

\Proof It suffices to prove that the sections of these vector
bundles coincide over the basis for the topology of $Z$ given by 
the $\phi ^{-1}(V)$ where $V$ runs over the affine  open subsets  of $Y$. 
Let $U=\phi ^{-1}(V)$ be one such.
Write   $R=O_Y(V)$,  $S=O_Z(U)$ and 
$M=E(V)\otimes _R\pi_*({\cal D}^{-1/2}_{X/Y})(V)$.
>From the very definitions we obtain that
$
\phi ^{*}(E_{\rho ,X})(U)=M^G\otimes _{R}S
$
and
$
\phi ^{*}(E)_{\phi ^{*}(\rho) ,T'}(U)=(\phi ^{*}(E)(U)
\otimes _{S}\pi^{(0)}_*({\cal D}^{-1/2}_{T'/Z})(U)))^G 
$.
We now use the fact that on the one hand
$
\phi ^{*}(E)(U)=E(V)\otimes _{R}S
$
while  on the other hand  by   [CNET1], Lemma~3.7 we know that
$
\pi^{(0)}_*({\cal D}^{-1/2}_{T'/Z})(U)=S\otimes _{R}
\pi_*({\cal D}^{-1/2}_{X/Y})(V)
$.
 Therefore we are reduced to proving that
$
(M\otimes _{R} S)^G=M^G\otimes _{R} S 
$
 and this follows at once since $S$ is flat over $R$.
\bigskip

\noindent To compare the Hasse-Witt invariants
of $\phi ^{*}(E)_{\phi ^{*}(\rho),T'}$ and $\phi ^{*}(E)_{\phi ^{*}(\rho),T}$ we now  use
the decomposition of the normalisation map
$$
T=T^{(m)}\rightarrow T^{(m-1)}\rightarrow \cdots \rightarrow
T^{(0)}=T' \ ,
$$
recalled in the previous sub-section,   
to construct a new family of metabolic  
bundles which allows us to use the Main Lemma. For any 
$0 \leq h \leq m$,  the group $G$ acts on $T^{(h)}$ and $\rho$ induces 
an orthogonal representation $\rho: G \rightarrow 
\Gamma (Z, {\bf O}(\phi^*(q))$. Therefore we can consider the twist of 
$(\phi^*(E),\phi^*(q))$  by this representation. For simplicity and when there is no ambiguity on the choice of the representation we will denote by  $ (\phi ^{*}(E)_{T^{(h)}},\phi ^{*}(q)_{T^{(h)}})$ the symmetric bundle 
$ (\phi ^{*}(E)_{\phi^*(\rho), T^{(h)}}, \phi^*(q)_{\phi^*(\rho), T^{(h)}})$. 
The principal advantage in considering the decomposition of normalisation into $m$ steps is that  we will be able
to compare 
the Hasse-Witt invariants of two consecutive terms $ \phi ^{*}(E)_{T^{(h)}}$ 
and $ \phi ^{*}(E)_{T^{(h+1)}}$. 

For $0\leq h\leq m-1$ we have the short exact 
sequence of locally free $O_Z$-modules
$$
0\rightarrow I^{(h)}\rightarrow \Lambda ^{(h)}\oplus \Lambda
^{(h+1)} \rightarrow {\cal G}^{(h)}\rightarrow 0 \ .
$$
Since $\phi ^{*}(E)$ is a locally free
$O_Z$-module  we obtain a new exact sequence
$$
0\rightarrow \phi ^{*}(E)\otimes _{O_Z} I^{(h)}\rightarrow \phi
^{*}(E)\otimes _{O_Z} \Lambda ^{(h)}\oplus \phi ^{*}(E)\otimes
_{O_Z}\Lambda ^{(h+1)} \rightarrow \phi ^{*}(E)\otimes _{O_Z}{\cal
G}^{(h)} \rightarrow 0 \ .
$$
>From the definitions of $I^{(h)}$ and  ${\cal G}^{(h)}$  it is
clear that these modules are $G$-modules and that the morphisms
in this sequence respect the action of $G$ when $G$ acts
diagonally on the tensor products. Next we consider the sequence
obtained  by taking  $G$-fixed points.

\begin{prop}
 For $0\leq
h\leq m-1$
$$
0\rightarrow (\phi ^{*}(E)\otimes _{O_Z} I^{(h)})^G \rightarrow
\phi ^{*}(E)_{T^{(h)}}\oplus \phi ^{*}(E)_{T^{(h+1)}} \rightarrow
(\phi ^{*}(E)\otimes _{O_Z}{\cal G}^{(h)})^G \rightarrow 0
$$
is an exact sequence of locally free $O_Z$-modules.
\end{prop}

\Proof It will be proved in Lemma 4.6, that $I^{(h)},\Lambda ^{(h)}$ 
and ${\cal G}^{(h)}$ are  locally projective $O_Z[G]$-modules for any 
$h$, $0\leq
h\leq m-1$. Therefore,  for any (closed) point $z$ of $Z$ we  
may choose a  sufficiently
small  affine neighbourhood $U$ of $z$ such that $I^{(h)}(U)$,
${\cal G}^{(h)}(U)$, $ \Lambda ^{(h)}(U)$ and $ \Lambda
^{(h+1)}(U)$ are  projective $O_Z(U)[G]$ -modules. 
>From  Lemma~\ref{ML} it follows   that 
the modules obtained by tensoring these modules with $\phi ^{*}(E)(U)$
over $O_Z(U)$  are all  projective 
$O_Z(U)[G]$-modules. Hence, by 
Lemma~\ref{fixedtrace}, we know that  
their $G$-fixed points  
are  locally free $O_Z(U)$-modules. We now consider the exact sequence 
of left $O_Z(U)[G]$-modules
$$0\rightarrow I^{(h)}(U)\rightarrow \Lambda ^{(h)}(U)\oplus 
\Lambda ^{(h+1)}(U)
\rightarrow {\cal G}^{(h)}(U)\rightarrow 0 \ .
$$
By tensoring with the $O_Z(U)$-locally free module $\phi ^{*}(E)(U)$
we obtain a new exact sequence of $O_Z(U)[G]$-modules, where $G$ acts
diagonally, 
$$
\begin{array}{lll}
0\rightarrow & \phi ^{*}(E)(U)\otimes _{O_Z(U)}I^{(h)}(U)  & \\
& & \\
 & \rightarrow
\phi ^{*}(E)(U)\otimes _{O_Z(U)}\Lambda ^{(h)}(U)\oplus \phi ^{*}(E)(U)
\otimes _{O_Z(U)}\Lambda ^{(h+1)}(U) & \\
& & \\
 & \ \ \ \ \ \ \ \
 \rightarrow \phi ^{*}(E)(U)\otimes _{O_Z(U)}{\cal G}^{(h)}(U)
 \rightarrow
0 \ .
\end{array}
  $$
 Taking $G$-fixed points affords a further   exact sequence
$$\begin{array}{lll}
0 \rightarrow (\phi ^{*}(E)(U)\otimes _{O_Z(U)} I^{(h)}(U))^G  \rightarrow
\phi ^{*}(E)_{T^{(h)}}(U)\oplus \phi ^{*}(E)_{T^{(h+1)}}(U)& \\
 & & \\
 \rightarrow
(\phi ^{*}(E)(U)\otimes _{O_Z(U)}{\cal G}^{(h)}(U))^G \rightarrow H^1(G, \phi ^{*}(E)(U)\otimes _{O_Z(U)} I^{(h)}(U)) \rightarrow 0 
 \ .
\end{array}
$$ 
Using Lemma~\ref{ML} and the fact  that $I^{(h)}(U)$ is 
a projective $O_Z(U)[G]$-module  we conclude that 
 $\phi ^{*}(E)(U)\otimes _{O_Z(U)} I^{(h)}(U)$ is also 
a projective  module over $O_Z(U)[G]$ and therefore that the last 
term of the previous sequence is equal to $0$. This last  exact sequence is
precisely the one we require  to complete  the proof of the proposition.
\medskip

\noindent {\bf Remark.} We observe that in the case where  
$E$ is a locally projective $O_Y[G]$-module it follows that 
$\phi ^{*}(E)$ is a projective $O_Z[G]$-module  and thus the proof  of 
the  proposition is complete without having to check the local 
projectivity of the modules $ I^{(h)},\Lambda ^{(h)}$ and ${\cal G}^{(h)}$. 

\bigskip

\noindent We deduce from Prop.~4.4  that for $0\leq h\leq m-1$ the
symmetric bundle obtained as the orthogonal sum of $(\phi
^{*}(E)_{T^{(h)}},(-1)^h\phi ^{*}(q)_{T^{(h)}})$ and \newline
$(\phi^{*}(E)_{T^{(h+1)}},(-1)^{h+1}\phi ^{*}(q)_{T^{(h+1)}})$ is metabolic 
with lagrangian \newline
$(\phi ^{*}(E) \otimes _{O_Z} {\cal
G}^{(h)})^G)^{\vee}$.
Hence, using the Main Lemma, we deduce that the product 
$$w_t(\phi
^{*}(E)_{T^{(h)}},(-1)^h\phi ^{*}(q)_{T^{(h)}}))w_t(\phi
^{*}(E)_{T^{(h+1)}},(-1)^{h+1}\phi ^{*}(q)_{T^{(h+1)}})$$
 is equal to $d_t((\phi ^{*}(E) \otimes _{O_Z} {\cal
G}^{(h)})^G) \ $. By taking the product of these equalities we then deduce that   $w_t(\phi
^{*}(E)_{T^{(0)}},\phi ^{*}(q)_{T^{(0)}}))$ equals to  
$$w_t(\phi
^{*}(E)_{T^{(m)}},(-1)^{m}\phi ^{*}(q)_{T^{(m)}})\prod _{0 \leq h
\leq m-1} d_t((\phi ^{*}(E) \otimes _{O_Z} {\cal
G}^{(h)})^G)^{(-1)^h} \ .$$
To deduce  Thm.~\ref{big-main-lemma} from this equality we simply
have to remember that $T^{(m)}=T$ by definition and that  $(\phi
^{*}(E)_{T^{(0)}},\phi ^{*}(q)_{T^{(0)}})=\phi^*(E_{\rho,X},q_{\rho,X})$  comes  from  Lemma 4.3.
\bigskip

\noindent We now establish the first part of Thm.~0.7 together with
a preliminary version of the second part. Then in the next section 
we shall conclude the proof of the second part. 
As a consequence of the fact that for $0\leq h\leq m$ the
restrictions of the forms $(\phi ^{*}(E)_{T^{(h)}},\phi
^{*}(q)_{T^{(h)}})$  to the generic fiber of $Z$  coincide,  we
deduce that  $w_1(\phi ^{*}(E)_{T^{(h)}})$ and  $c_1(\phi
^{*}(E)_{T^{(h)}})$ do not depend on the choice of $h$. 
Moreover it also follows from  Lemma 4.2 $(c)$ that
$$
w_2(\phi ^{*}(E)_{T^{h}})-w_2(\phi ^{*}(E)_{T^{h+1}})
=c_1((\phi ^{*}(E)\otimes _{O_Z}{\cal G}^{(h)})^G )-c_1(\phi ^{*}(E)_{T^{h+1}})
$$
 for $0\leq h\leq m-1$. Therefore by adding these equalities we obtain that
$$ 
w_2(\phi ^{*}(E_{\rho,X}))-w_2(\phi ^{*}(E)_{T})= \delta(\phi ^{*}(E),Z)
$$
 with
 $$
\delta(\phi ^{*}(E),Z)=
\sum _{ 1\leq h\leq m}\left[
c_1((\phi ^{*}(E)\otimes _{O_Z}{\cal G}^{(h-1)})^G )-c_1(\phi ^{*}(E)_{T^{(h)}})\right]\ .
$$
Using Thm.~\ref{etale} to evaluate 
$w_k(\phi ^{*}(E)_{T})-w_k(\phi ^{*}(E))$ 
we finally obtain that
$$
w_1(\phi ^{*}(E_{\rho,X}))=w_1(\phi^*(E)_T)=w_1( \phi ^{*}(E)) +w_1(\phi ^{*}(\rho))
$$
and
\begin{eqnarray*}
w_2(\phi ^{*}(E_{\rho,X})) & =& w_2( \phi ^{*}(E))+
w_1( \phi ^{*}(E))w_1( \phi ^{*}(\rho))+ w_2(\phi ^{*}(\rho))+\\
 \ \ \ \ \ & + & sp_2(\phi ^{*}(\rho))+ \delta(\phi ^{*}(E),Z) \ .
\end{eqnarray*}

\bigskip

\noindent {\bf Remark.} Since $2$ is invertible in $Z$, we may
associate to the squaring map on ${\bf G}_m $ an exact
 Kummer  sequence of \'etale sheaves of groups
$$ 
0\rightarrow {\bf Z}/2{\bf Z} \rightarrow {\bf G}_m
\rightarrow {\bf G}_m \rightarrow 0 
$$ 
and therefore an exact
sequence of groups
$$  
0\rightarrow \mbox{\rm Pic}(Z)/2 \rightarrow H^2(Z,{\bf Z}/2{\bf  Z})
 \rightarrow   H^2(Z,{\bf G}_m)_2 \rightarrow 0 \ .
$$
The group $H^2(Z,{\bf G}_m)$ is known as the cohomological Brauer
group of $Z$ and is denoted by $Br'(Z)$, (see for instance 
[Mi], Chapter $4$). It  now follows from   the  definition  of the first Chern
class that $ \delta(\phi ^{*}(E),Z)$ belongs to $\mbox{\rm Pic}(Z)/2$ 
(see below). 
Therefore
by projecting our formula into  $Br'(Z)$ we obtain a   formula of
Fr\"ohlich-type ([F], Thms 2 and 3),  in this group. 
It should be  observed that in our
case--$Z$ regular and integral--one believes that $Br'(Z)$  and
the Brauer group of $Z$ coincide.  Furthermore, since   $Z$ is
integral and regular, we have an exact sequence, \cite{EKV} 5.4
$$  
0\rightarrow \mbox{\rm Pic}(Z)/2 \rightarrow H^2(Z,{\bf Z}/2{\bf  Z})
 \rightarrow  H^2(K,{\bf Z}/2{\bf Z}) \rightarrow 0 \ ,
$$
where $K$ denotes  the function field of $Z$. Therefore,  when we restrict
our formul{\ae} to the generic fiber of $Z$,  the term $
\delta(\phi ^{*}(E),Z)$ disappears and so we obtain Fr\"ohlich's
original formula in the Galois cohomology  group $H^2(K,{\bf
Z}/2{\bf Z})$.
\bigskip

\subsection{Proof of Theorem~\ref{tame}}

\noindent To complete the proof of Thm.~\ref{tame} we shall   
often  need to refer to
[CNET1], Sect.~$4$.  We first observe that, from  the definition
of the first   Chern class, it follows  that  for any $ 1\leq h\leq m
$, the element of $H^2(Z, {\bf Z}/2{\bf Z})$
$$  c_1((\phi ^{*}(E)\otimes _{O_Z}{\cal G}^{(h-1)})^G )-
c_1(\phi ^{*}(E)_{T^{(h)}}) $$  is the image  of the element $
[ \mbox{\rm det}(\phi ^{*}(E)\otimes _{O_Z}{\cal G}^{(h-1)})^G )]-
[\mbox{\rm det}( \phi ^{*}(E)_{T^{(h)}})]  $ of $\mbox{\rm Pic}(Z)$, where we
denote by $[D]$  the class in $\mbox{\rm Pic}(Z)$ of the divisor $D$ of $Z$.
Let $\alpha^{(h)} :   \phi ^{*}(E)_{T^{h}} \rightarrow  (\phi
^{*}(E)\otimes _{O_Z}{\cal G}^{(h-1)})^G$ be the inclusion map
 (see Prop.~4.4). We
have the short  exact sequence of $O_Z$-modules
$$
0 \rightarrow \mbox{\rm det}( \phi ^{*}(E)_{T^{(h)}})
\stackrel{\mbox{\scriptsize  det}(\alpha ^{(h)})}{\rightarrow} 
\mbox{\rm det}\left( (\phi ^{*}(E)\otimes _{O_Z}{\cal
G}^{(h-1)})^G\right) 
\rightarrow \mbox{\rm coker} (\mbox{\rm det}(\alpha
^{(h)})) \rightarrow 0\  .
$$ 
Therefore there exists a divisor $\Delta ^{(h)}(E)$ such that
$$  
[\Delta ^{(h)}(E)]= [\mbox{\rm det}( \phi ^{*}(E)
\otimes _{O_Z}{\cal G}^{(h-1)})^G]
 -[\mbox{\rm det}( \phi ^{*}(E)_{T^{(h)}})]  \ .
$$
 One of the main ingredients  of the proof of the theorem is to  
provide  an explicit description of  each $\Delta ^{(h)}(E)$  
and thus of the sum  $\Delta (E)= \sum _{h=1}^{m}\Delta ^{(h)}(E)$  
at least when restricted to each $U_i$ of some \'etale 
covering $(U_i \rightarrow Z)$ of $Z$.  Since  $\Delta (E)$ 
has  image $\delta(\phi ^{*}(E),Z)$ in $H^2(Z, {\bf Z}/2{\bf Z})$,  
the theorem will follow from the congruence in $\mbox{\rm Div}(Z)$
$$ 
\Delta (E) \equiv \phi^*(R(\rho,X))\ \ mod\ 2 \ .
$$
This congruence will be deduced from the congruence mod $2$ of the
restrictions  of both terms to each $U_i$ of the covering.  The
final part of  this section is devoted to the construction for any (closed) point
of $Z$ of an \'etale neighbourhood on which we can evaluate and
compare  as divisors the restrictions of $\Delta (E)$ and $
\phi^*(R(\rho, X))$.

Before proceeding to these constructions and computations,  we start
by   fixing  some notation and by  recalling  the
results of [CNET1] that we shall require. Let $y$ be a 
point of $Y$ of arbitrary dimension and $\sigma :S \rightarrow Y$ be an \'etale
neighbourhood of $y$. For any scheme $v:V \rightarrow Y$ we  write
$V_S=V \times _Y S$, we   define $O_V(S)$ as the ring of its
global sections and   denote by $v_S : V_S \rightarrow S$ and
$\sigma _V :V_S \rightarrow V$  the projection maps. For a  $V$-vector bundle ${\cal
F}$ we obtain the vector bundle ${\cal F}
_S ={\sigma _V}^*({\cal F})$ over $ V_S$. We denote by  ${\cal
F}(S)$   the module of its  global sections.  Let now $z$ be point of $Z$ and $y=\phi (z)$. The \'etale neighbourhoods of
$z$ that we shall consider  will always be of the form $\sigma _Z
:Z_S \rightarrow Z$ where $\sigma :S \rightarrow Y$ is  a well-chosen  
\'etale  affine neighbourhood of $y$. As  mentioned
previously,   the objects that we  consider all have good functorial
properties under  base change. To be  more precise:  suppose that
$\sigma _Z :Z_S \rightarrow Z$ is an \'etale neighbourhood of the
type introduced above. For any $0 \leq h \leq m$ we can
consider  on the one hand $(T^{(h)})_S$ as defined previously  and
on the other hand  $T_S^{(h)}$ as the normalisation of $T'_S$ along
the divisor $\sigma^* (b_1)\cup ...\cup \sigma^* (b_h)$.  From
the fact that they coincide ([CNET1], Prop.~3.6), 
 we deduce as in Thm.~\ref{itsabundle} that
$$ ( \phi ^*(E)_{T^{(h)}} )_S=(( \phi ^*(E))_S)_{T_S^{(h)}}  \ , $$
and that we have an exact sequence of $Z_S$-vector bundles
$$0\rightarrow (F\otimes _{O_{Z_S}} I_S^{(h)})^G\rightarrow F_{T_S^{(h)}}\oplus   F _{T_S^{(h+1)}}
\rightarrow (F\otimes _{O_{Z_S}}{\cal G}_S^{(h)})^G \rightarrow 0
\ ,
$$
where for simplicity we denote by $F$ the $O_{Z_S}$-vector bundle
$\phi ^*(E)_S$. Therefore, using the fact that base-change
commutes with taking determinants, for any $1\leq h \leq m$, 
the restriction  ${\sigma _Z}^* (\Delta ^{(h)}(E))$
of $\Delta ^{(h)}(E)$ to the \'etale neigbourhood
$Z_S$ of $z$
will be obtained via the exact sequence
$$0 \rightarrow \mbox{\rm det}( F_{T_S^{(h)}})
\rightarrow \mbox{\rm det}( F\otimes _{O_{Z_S}}{\cal
G}_S^{(h-1)})^G \rightarrow \mbox{\rm coker} (\mbox{\rm
det}(\alpha_S^{(h)})) \rightarrow 0\  ,$$ with   $\alpha_S^{(h)} :
F_{T_S^{(h)}} \rightarrow  (F\otimes _{O_{Z_S}}{\cal
G}_S^{(h-1)})^G$ again being the inclusion map.

We now want to make precise the choice of  the neigbourhood $S$
 which allows us to obtain explicit descriptions of the global sections of  
$ F_{T_S^{(h)}}$ and

\noindent $ (F\otimes _{O_{Z}(S)}{\cal G}_S^{(h-1)})^G$ which lead
themselves
 to determinantal computations.
\bigskip

We fix $z$ a  point of $Z$, we choose a point $x$ of $X$
whose image in $Z$ is $z$ and we write $y$ for $\phi (z)$.  For any 
$1\leq h \leq m$,  the inertia group of a generic point of $
B_{h,k}$ only depends on $h$ up to conjugacy, thus  we may  denote
by $I_h$ this group and by $e_h$ the  order of $I_h$.  The cover  $X/Y$
being  tame,  $I_h$ is cyclic and  therefore will often be
identified with ${\bf Z}/e_{h}{\bf Z}$. 
Let $I({ x})$ be the  inertia group of  the point $x$.  Then we
know (see for instance Sect.~$2$ of [CNET1]) , that
$$ 
I({ x}) \cong  \prod _{\ell \in J({ x})}I_{\ell}
\cong \prod _{\ell \in J({ x})}{\bf Z}/e_{\ell}{\bf Z} \ , 
$$
where
$$
J({ x}) =  \{ \ell \ | \ 1\leq \ell \leq m, \exists k : { x} \in
B_{\ell,k} \}   \ .
$$
After reordering if necessary we shall take $J({ x}) =  \{ 1, 2,
\ldots , n \}$. For any such $\ell$, we  denote by $\chi
_{\ell}$ the character giving the action of $I_{\ell}$ on the
cotangent space at the generic point of $B_{\ell,k}$. It follows
from  [CNET1], Lemma 2.3, that there exists an integral affine
\'etale neighbourhood of $y$,
$$
S=\mbox{\rm Spec}(A_y) \ ,
$$ 
where $A_y$ is an algebra containing a sequence $a_1, a_2, \ldots,
a_n$ of regular parameters and enough roots of unity of order
coprime with the residue characteristic of $y$ and an isomorphism
of schemes with $G$-action
$$
X_S \cong \mbox{\rm Spec}(O_{ X}(S))
$$
where
$$
O_{ X}(S) :=\mbox{Map}_{I({ x})}(G,B_{ x})
$$ and
$$
B_{ x}:= A_y[t_1,\ldots ,t_n] = A_y[T_1,\ldots
,T_n]/(T_{1}^{e_1}-a_1,\ldots , T_{n}^{e_n}-a_n) \ .
$$ 
Moreover, for $ 1\leq \ell \leq n$, the action of $I_{\ell}$  
on  the image of $T_l$ in $B_{x}$ that we denote by $t_{\ell}$,   
is given by the character
  $\chi _{\ell}$. We may now obtain  a  description of $Z_S$ and a description  of   
$T'_S$  considered as a scheme with $G$-action,  from the above description  of $X_S$.
More precisely, if $s$ denotes the cardinality of $G_{2}\backslash
G/I(x)$,  we obtain that
\begin{eq} \label{ZS}
Z_S = \mbox{Spec}(O_Z(S)) \ \ \ \mbox{with} \ \ \ O_Z(S)=\prod
_{1\leq j\leq s}B_{j,x} \ ,
\end{eq}
which is a product of $s$ copies of $B_{j,x}=B_{x}$
 and
$$
T'_S = \mbox{Spec}(O_{T'}(S)) \ \ \ \mbox{with} \ \ \
O_{T'}(S)=\mbox{Map}_{I({ x})}(G, \prod _{ 1\leq j\leq s}
B_{j,x}\otimes _{A_y}B_{x}) \ ,
$$ and $I(x)$ acting on the second factor of the tensor products.  From now on the schemes  $S$ and
$Z_S$ given above will be the neighbourhoods of $y$ and
$z$ that we shall consider.

 We denote by ${\cal {S}} (x) $ the set of
sequences $\alpha=(\alpha _\ell),\  \ell \in J({ x}) $ where    for
each $\ell \in J({ x}) $ we have an integer  $\alpha _\ell $ such that $0\leq \alpha _\ell< e_{\ell}$. For any $1 \leq h
\leq m$, we consider the  partition of $J=J(x)$ into
$$
J'_h = J'_h(x) = \{ \ell \in J : 1\leq \ell\leq h \}
$$
and
$$
J''_h = J''_h( x) = \{ \ell \in J : h+1\leq \ell \leq m \} \ .
$$
So for $h\geq n$, the set $ J''_h$ is empty. For $\alpha = (\alpha
_{\ell})$, we  write
$$
\partial ^{(h)}(\alpha ) := \left(
\prod _{\ell\in J'_h} t^{-\alpha _{\ell}}_{\ell}
             \prod _{\ell\in J''_h} t^{-\epsilon _{\ell}}_{\ell} \right)
                  \otimes \prod _{\ell\in J} t^{\alpha _{\ell}}_{\ell} \ ,
$$
where $\epsilon _{\ell} = 0$ or $e_{\ell}$ depending on whether
$\alpha _{\ell}$ is strictly smaller or strictly larger than
$e_{\ell}/2$. Then we define 
$$
D^{(h)}_{j}(\alpha ) = B_{j,x}\partial ^{(h)}(\alpha ) \ \ \
\mbox{and} \ \ \ D^{(h)}_{j}(S) =
                \bigoplus _{\alpha \in {\cal {S}} (x)} D^{(h)}_{j}(\alpha )  
\ ,
$$ 
and
$$
{\cal G} ^{(h-1)}_{j}(\alpha)=D^{(h-1)}_{j}(\alpha )+
D^{(h)}_{j}(\alpha )\ \ \mbox{and}\ \ {\cal G} ^{(h-1)}_{j}(S) =
\bigoplus _{\alpha \in {\cal {S}} (x)} {\cal G}
^{(h-1)}_{j}(\alpha ) \ .
$$
Let us denote by $\chi _{x}^{\alpha}$ the character $\prod _{ \ell
\in J(x)} \chi _{\ell}^{\alpha_\ell}$ of $I(x)$;  then we shall
observe that the previous decompositions correspond to the
decomposition according to the characters of $I(x)$. From now on,
for any $I(x)$-module $M$,  we denote by $M(\alpha)$ the submodule
of $M$ on which the action of $I(x)$ is given by $\chi
_{x}^{\alpha}$.  It now follows from [CNET1], Sect.~3.e. that we
have the isomorphisms of $G$-modules
$$
I ^{(h)} (S) \cong \mbox{Map}_{I({ x})}(G, \prod _{ 1\leq
j\leq s} I^{(h)}_{j}(S))
 \ ,$$
$$\Lambda ^{(h)} (S) \cong \mbox{Map}_{I({ x})}(G, \prod _{ 1\leq
j\leq s} D^{(h)}_{j}(S))
  $$
and
$$
{\cal G}
^{(h-1)}(S)\cong \mbox{Map}_{I({ x})}(G, \prod _{ 1\leq j\leq s}
{\cal G}^{(h-1)}_{j}(S)) \ .
$$ 
We start by deducing from  these descriptions, 
that $I^{(h)}$, $\Lambda ^{(h)}$ and ${ \cal G}
^{(h)}$ are locally projective $O_Z[G]$-modules
(recall that this was used in the proof of Prop.~4.4).

\begin{lemma} 
For any $h$, $1\leq h \leq m$, the $O_Z[G]$-modules  
$I^{(h)}$, $\Lambda ^{(h)}$ and ${ \cal G}^{(h)}$ are locally projective.
\end{lemma}
\Proof Let $z$ be a  point of $Z$, $y=\phi (z)$ and  
let $\sigma: S \rightarrow Y$ be the \'etale neighbourhood of $y$ 
introduced above. Since $\sigma$ is of finite type it is open. 
Let us denote by $U$ the image $\sigma(S)$. We may 
assume that $U$ is affine. In fact, since $ \sigma $ is \'etale and of 
finite type, it is smooth and quasi-finite, [Mi], I, Remark 3.25. 
Therefore it follows from Zariski's Main Theorem that it can be decomposed 
as the product of an open immersion and a finite map. Hence we are reduced
 to the case where $\sigma  :S \rightarrow U$ is finite.  It now 
follows from a theorem of Chevalley \cite{EGAII}, Thm.~6.7-1, 
that $U$ is affine. 
Let us denote by $V$ the affine open neighbourhood $\sigma^{-1}(U)$ of 
$z$ and let us prove for instance  that $I^{(h)}(V)$ is a projective 
$O_Z(V)[G]$-module (the proof of the projectivity of 
$\Lambda ^{(h)}(V)$ and ${ \cal G}^{(h)}(V)$ will follow exactly 
the same lines). Let  $\sigma_Z : Z_S \rightarrow Z$ be the morphism 
obtained from $\sigma $ by base change. It is \'etale and moreover 
$\sigma_Z(Z_S)=V$. It follows that $\sigma_Z$ induces a faithfully 
flat morphism of affine schemes from $Z_S$ onto $V$. Therefore it 
suffices to prove that $I^{(h)}(S)$ is a projective $O_Z(S)[G]$-module. 
>From  the previous description of $I^{(h)}(S)$ we observe that this 
module is induced from the $O_Z(S)[I(x)]$-module  $\prod _{ 1\leq
j\leq s} I^{(h)}_{j}(S)$. From [CNET1], section 3.e., we know 
that this last module is finitely generated and free as an $O_Z(S)$-module. 
Since the order of $I(x)$ is a unit in $O_Z(S)$,  this is enough to conclude  
that it is projective as an  $O_Z(S)[I(x)]$-module and hence that 
$I^{(h)}(S)$ is a projective $O_Z(S)[G]$-module. 
This completes the proof of the lemma. 

\bigskip 
 We now return to the proof of Thm.~0.7. 
We recall that $F$ has been defined at the beginning of this section as 
the $O_{Z_S}$-vector bundle $\phi(E)_S$. Denoting respectively by $F(S)$,
$E(S)$  and   $F_{T_S^{(h)}}(S)$ the global sections of $F$, 
$\varphi ^*(E)$ and  $F_{T_S^{(h)}}$,  from the previous
equalities we deduce that
 $$ 
(F(S)\otimes _{O_Z(S)}{\cal G}^{(h-1)}(S))^G \cong
\left( \mbox{Map}_{I({ x})}(G, \prod _{ 1\leq j\leq s}
 (E(S)\otimes _{A_y}{\cal G}^{(h-1)}_{j}(S)) \right) ^G 
$$
and
$$ 
F_{T_S^{(h)}}(S)\cong 
\left( \mbox{Map}_{I({ x})}(G, \prod _{ 1\leq j\leq s}
(E(S)\otimes _{A_y}D^{(h)}_{j}(S))\right) ^G \ .$$ For any $R[G]$-module $M$ and any subgroup $H$ of $G$ we recall that we denote by  $\mbox{Map}_{H}(G,M)$ the set of maps 
$u:G \rightarrow M$ such that for any $x\in G$ and any $h \in H $, $u(hx)=hu(x)$,  endowed with its structure of $G$-module defined by 
$gu:x \mapsto u(xg) \ .$ It is now easy to
check that for any such module $M$ and any  such subgroup $H$ of $G$ the
map $f \mapsto f(1)$ induces  an isomorphism of $R$-modules
from $(\mbox{Map}_{H}(G,M))^G$ onto $M^H$. Therefore we conclude
that we have the following isomorphisms
$$
(F(S)\otimes _{O_Z(S)}{\cal G}^{(h-1)}(S))^G \cong \prod _{ 1\leq
j\leq s} (E(S)\otimes _{A_y}{\cal G}^{(h-1)}_{j}(S)))^{I({ x})}$$
and
$$
F_{T_S^{(h)}}(S)\cong \prod _{ 1\leq j\leq s} (E(S)\otimes
_{A_y}D^{(h)}_{j}(S)))^{I(x)}\ .
$$
Denoting by  $E_j(S)$ the tensor product $E(S) \otimes _{A_y}
B_{j,x}$, we finally obtain that
$$
(F(S)\otimes _{O_Z(S)}{\cal G}^{(h-1)}(S))^G \cong \prod _{ 1\leq
j\leq s} (E_j(S)\otimes _{B_{j,x}}{\cal G}^{(h-1)}_{j}(S)))^{I({
x})}$$ and
$$
F_{T_S^{(h)}}(S)\cong \prod _{ 1\leq j\leq s} (E_j(S)\otimes
_{B_{j,x}}D^{(h)}_{j}(S)))^{I(x)}\ .
$$
We now want to consider the structure of $E$ as a module when
restricted to $I(x)$. To that end  we consider  $E_y$ as an  $O_{Y,y}[I(x)]$-module. Since the order of $I(x)$ and the residue characteristic of $y$ are coprime, $E_y$ can be decomposed according to the characters of $I(x)$ as a direct sum
$$
E_y= \sum _{\alpha \in {\cal S}(x)} E_y(\alpha) \ .
$$ 
For any $\alpha$,  let us denote by $l_y(\alpha)$ the  
$O_{Y,y}$-rank  of $E_y(\alpha)$. 
Since $E_y$ is an orthogonal representation of $I(x)$
we observe that   for any
$\alpha\in {\cal S} (x)$ it follows that  $ l_y(\alpha)=l_y(e-\alpha)$ where
$e=(e_{\ell}), \ \ \ell \in J(x)$. Therefore for any $1 \leq
j \leq s$, we obtain  a decomposition for $E_j(S)$, namely
$$E_j(S)= \sum _{\alpha \in {\cal S}(x)} E_j(S)(\alpha) \ ,$$ where $E_j(S)(\alpha)$ is a free $B_{j,x}$-module of  rank $l_y(\alpha)$. Hence we deduce that
$$(E_j(S)\otimes_{B_{j,x}} {\cal G}^{(h-1)}_{j}(S)))^{I({ x})} \cong
\oplus _{\alpha \in {\cal S}(x)} (E_j(S)(\alpha)\otimes
_{B_{j,x}}{\cal G}^{(h-1)}_{j}(S)(e-\alpha)) $$ and
$$(E_j(S)\otimes_{B_{j,x}} D^{(h)}_{j}(S)))^{I({ x})} \cong
\oplus _{\alpha \in {\cal S}(x)} (E_j(S)(\alpha)\otimes
_{B_{j,x}}D^{(h)}_{j}(S)(e-\alpha))\ .$$

We now return  to the computations of determinants. For $h
\notin {J}(x) $ and  for any $j$ and $\alpha$ we know from [CNET1], 
Prop.~3.14, that $
D^{(h)}_{j}(S)(\alpha)=D^{(h-1)}_{j}(S)(\alpha)$ and so 
is equal to ${\cal G}^{(h-1)}_{j}(S)(\alpha)$. Therefore  
$$ 
\mbox{\rm det}( F_{T_S^{(h)}}) = \mbox{\rm det}(
F\otimes _{O_{Z_S}}{\cal G}_S^{(h-1)})^G \ .
$$ 
 Assuming  now that $h
\in {J}(x)$  we introduce the partition of $ {\cal S}(x)$ into
${\cal S}_h(x)$ and ${\cal S}'_h(x)$ where ${\cal S}_h(x)$ 
(resp.  ${\cal S}'_h(x)$) denotes
the set of sequences $\alpha$ such that $e_h/2 <\alpha_h$ (resp. $\alpha_h<e_h/2$). We
recall from [CNET1], Prop.~3.14 that
$D^{(h)}_{j}(S)(\alpha)=(t_h^{(e_h-\alpha_h)}\otimes 1)\  {\cal
G}^{(h-1)}_{j}(S)(\alpha)$ for $\alpha \in  {\cal S}_h(x)$ and
equal to ${\cal G}^{(h-1)}_{j}(S)(\alpha)$ otherwise. 
Hence  we deduce from these equalities that, if 
$\alpha \in {\cal S}_h(x)$, then $(e-\alpha) \in {\cal S}'_h(x)$ 
and thus we have   
$$
\mbox{\rm det}(E_j(\alpha)\otimes _{B_{j,x}}{\cal G}^{(h-1)}_{j}(S)(e-\alpha))
=\mbox{\rm det}(E_j(\alpha)\otimes _{B_{j,x}}D^{(h)}_{j}(S)(e-\alpha)) \ . 
$$
 If now  $\alpha \in {\cal S}_h'(x)$, then   
$(e-\alpha) \in {\cal S}_h(x)$ and 
therefore  $$\mbox{\rm det}(E_j(\alpha)\otimes _{B_{j,x}}D^{(h)}_{j}(S)(e-\alpha))=
\mbox{\rm det}(E_j(\alpha)\otimes _{B_{j,x}} t_h^{\alpha _h}{\cal G}^{(h-1)}_{j}(S)(e-\alpha)) \ .
$$
Using  the fact that the
$E_j(S)(\alpha)$ are free $B_{j,x}$-modules  of  rank $l_y(\alpha)$ and that  $D^{(h)}_{j}(S)(\alpha)$  and  ${\cal
G}^{(h-1)}_{j}(S)(\alpha)$ are free $B_{j,x}$ rank one modules
([CNET1], Prop.~3.14),   we deduce
from  above  that in this last case  
$$
\mbox{\rm det}(E_j(\alpha)\otimes _{B_{j,x}}{D}^{(h)}_{j}(S)(e-\alpha))=
\mbox{\rm det}(E_j(\alpha)\otimes _{B_{j,x}}{\cal G}^{(h-1)}_{j}(S)(e-\alpha))t_h^{\alpha _h l_y(\alpha)} \ .
$$
Thus for  any $ 1 \leq j \leq s$,
$$
\mbox{\rm det}((E_j(S)\otimes _{B_{j,x}}D^{(h)}_{j}(S))^{I(x)}) =
 \mbox{\rm det}( (E_j(S)\otimes _{B_{j,x}}{\cal G}^{(h-1)}_{j}(S))^{I(x)})\ t_h^{\sum _{\alpha \in {\cal S}'_h(x)}
\alpha_h l_y(\alpha)} \ ,
$$
 and therefore from the previous equalities  we  may finally conclude
 that 
$$  
\mbox{\rm det}( F_{T_S^{(h)}}) = \mbox{\rm det}(( F\otimes _{O_{Z_S}}{\cal G}_S^{(h-1)})^G)
\prod _{ 1\leq j\leq s} t_h^{\sum _{\alpha \in {\cal S}'_h(x)}
\alpha_h l_y(\alpha)} \ . 
$$ 
It follows from this equality that,  
for $h \in J(x)$,  the restriction of $\Delta ^{(h)}(E)$ to $Z_S$ 
is defined as a Cartier divisor by
$$ \Gamma ^{(h)}(E)= \prod _{ 1\leq j\leq s} t_h^{-\gamma^{(h)}(E)} \ ,$$ 
where we write $\gamma^{(h)}(E)={\sum _{\alpha \in {\cal S}'_h(x)}\alpha_h l_y(\alpha)}$.

 We now want to give an interpretation of $ \gamma^{(h)}(E)$. We start by observing that,  from its very definition, 
$$
\gamma^{(h)}(E)= \sum _{0 \leq k \leq e_h/2}k 
\sum _{\alpha \in {\cal S}_h(x)\atop  \alpha_h=k} l_y(\alpha) \ .
$$
For any  $1 \leq h\leq m$,  let us  denote by  $d_k^{({h})}(E)$ the rank over $O_{Y,\xi _{h}}$ of the $\chi_{h}^k$ component of $E_{\xi _{h}}$, when considered as an $I_h$-module. When $h \in J(x)$,  since $O_{Y,\xi _{h}}$ contains $O_{Y,y}$,  we can write  $E_{\xi _{h}}$ as the tensor product $E_y\otimes _{O_{Y,y}}O_{Y,\xi _{h}}$. We  then deduce from this equality the following decomposition of  $E_{\xi _{h}}$ into a direct sum of $I(x)$-modules:  
$$ E_{\xi _{h}}= \oplus _{\alpha \in {\cal S}(x)} ( (E_y)(\alpha)\otimes_{O_{Y,y}}O_{Y,\xi _{h}}) \ .$$
This therefore  implies that the $\chi_{h}^k$-component of $E_{\xi _{h}}$, when considered as an $I_h$-module,  is the direct sum of the $E_y(\alpha)\otimes_{O_{Y,y}}O_{Y,\xi _{h}}$ when $\alpha$ runs through the elements of $ {\cal S}_h(x)$ such that $\alpha_h=k$.  It then follows from this  decomposition that, for $ 0 \leq k < e_h/2$,  
$$
d_k^{(h)}(E)= \sum _{\alpha \in {\cal S}_h(x)\atop  \alpha_h=k} l_y(\alpha) 
\ .
$$
We have then proved that for $h \in J(x)$
$$
\gamma^{(h)}(E)=d^{({h})}(E) \ ,
$$
where $d^{({h})}(E)$ has been defined in the introduction 
as the sum $\sum _{k=0}^{e_{h}/2}k d_k^{({h})}(E)$.

We now come back to the function $\Gamma^{(h)}(E)$ which defines $\Delta^{(h)}(E)$. Writing $N = \prod _{h \in { J}(x)} e_h$ and, for $h \in
{ J}(x)$, $N_{h}=N/e_{h}$ and using that $t_h^{e_h}=a_h$ for any such $h$,  we obtain that 
$$
\Gamma ^{(h)}(E)^N =  a_{h}^{-N_{h}d^{(h)}(E)}\ .
$$
Since by hypothesis the ramification indices  are odd we conclude that 
the restriction of
$\Delta (E)$ to $Z_S$, namely $\sigma_Z ^*(\Delta (E))$,  is
defined as a Cartier divisor by the function
 $$
\Gamma (E)=  \prod _{h  \in { J}(x)}\Gamma^{(h)}(E) 
\equiv  \prod _{h \in { J}(x)} a_{h}^{d^{(h)}(E)}
\ \ \bmod 2\ .
$$ 
Let us now consider the ramification divisor
$$
R(\rho, X)=\sum _{1 \leq h \leq m} d^{(h)}(E)b_h  
$$
defined in the introduction. We write $U$ for the image of $\sigma $. We assume that $S$
has been chosen sufficiently small  such that the irreducible components
that $U$ intersects are precisely those containing  $y$, namely
$\{ b_h, h\in { J}(x) \}$. Since $\sigma : S   \rightarrow
U$ is \'etale,  it follows that
$$ 
\sigma ^*(R(\rho,X))= \sum _{h \in {\cal J}(x)} d^{{h }}(E)\sum _{\varphi (\eta)=
\xi_{h}} \{ \bar {\eta} \} \  ,
$$
where the $\eta$'s are  points of $S$ of codimension $1$ over $\xi_h$. 
Therefore,  since  $a_{h}$ is a local equation of $b_{h}$
for any such $h$, we obtain  that $ \sigma ^*(R(\rho,X))$ and
therefore $\phi_S ^*(\sigma ^*(R(\rho,X))$ is defined  by the
function $\prod _{h \in { J}(x)} a_{h}^{d^{(h)}(E)}$.
Using the equality $\phi_S ^*(\sigma ^*(R(\rho,X)(E))=
\sigma_Z^*(\phi ^*(R(\rho,X))$ and the congruence satisfied by
the function $\Gamma (E)$ we conclude that the restrictions of
$\phi ^*(R(\rho,X))$  and $\Delta (E)$ to $Z_S$ are indeed congruent mod $2$ as required.

\bigskip

\Example  
Our final goal is  to compute  the divisor $R(\rho,X)$ 
in a special case considered in Sect.~\ref{subquotients}. 
We keep the hypotheses and the notations of Example 2.6.  
So we consider the symmetric bundle $(E,q)$  
where $E=O_Y[G/H]$ and $q$ is the symmetric form which has 
the cosets $\{{\bar a}=aH \}$ as an orthonormal basis; 
 $\rho$ is the tame orthogonal representation of $G$ induced 
by permuting the ${\bar a}$. 

By the above work  the computation of this divisor $R(\rho,X)$ reduces to  evaluating the integers $d_k^{(h)}(E)$ for any 
$1 \leq h \leq m$ and $1 \leq k \leq e_h/2$. We now fix such an $h$,
 we  choose once for all a codimension one point $\xi_h''$ of $X$
above $\xi_h$ and we   assume for simplicity that $O_{Y,\xi_h}$ 
contains the values of the character $\chi_h$. We let $I_h$ 
(resp. $\Delta_h$) denote the inertia group (resp. decomposition group)
of $\xi_h''$. Since there is no risk of ambiguity we make no further  
mention from now of the dependance upon $h$ of the objects we consider 
and therefore we will write $I$ for $I_h$, $\Delta$ for $\Delta_h$, 
$\xi$ for $\xi_h$,  {\it etc.}

If $G=\cup_{1\leq i\leq r} \Delta \gamma_iH$ is a double coset 
 decomposition of $G$ then, by standard theory (see for instance 
[FT], Chap.~8, Sect.~7) we have an isomorphism of 
left $O_{Y,\xi}[\Delta]$-modules
$$
E_{\xi} \cong \oplus_{ 1\leq i\leq r} O_{Y,\xi}[\Delta /(^{\gamma_i}H\cap \Delta)] \ ,
$$ 
where for $\gamma \in G$ we write $^{\gamma}H=\gamma H
{\gamma}^{-1}$. We observe that the above double cosets parametrise
the codimension one points of $V=X/H$ above $\xi$ by the rule
$\gamma_i \mapsto \lambda (\xi^{\gamma_i})$ (we recall that $G$
acts on the right on $X$). Then $\Delta_i=H\cap
^{\gamma_i^{-1}}\Delta$ (resp. $I_i=H\cap ^{\gamma_i^{-1}} I)$ is the
decomposition group (resp. the inertia group) of $\xi''^{\gamma_i}$
over $V$. We recall that we write $e$ (resp. $f$) for the
ramification index (resp. the residue class extension degree)  of
$\xi''$ over $Y$. Let us write $e'_i$ (resp. $f'_i$) for the
ramification index (resp. the residue class extension degree)  of
$\xi''^{\gamma_i}$ over $V$, thus the codimension one point on $V$
corresponding to $\gamma_i$, namely $\lambda(\xi^{\gamma_i})$,  has
ramification $e_i=e{e'_i}^{-1}$ and residue class extension degree
$f_i=f{f'_i}^{-1}$. It follows from the above isomorphism that in
order to decompose $E_{\xi}$ as a direct sum of $O_{Y,\xi}[I]$-modules
we have to decompose each $O_{Y,\xi}[\Delta /((^{\gamma_i}H\cap
\Delta)]$. With this in mind we consider the double cosets  
decomposition of 
$I\backslash \Delta / (^{\gamma_i}H\cap \Delta)$. Using the fact
 that $I$ is a normal subgroup of $\Delta$ we observe that each component of the direct sum decomposition is isomorphic to 
$O_{Y,\xi}[I/(^{\gamma_i}H\cap I)]$ 
and that the number of components, 
namely the number of the sets of double cosets,  is equal to $f_i$. 
Therefore we have proved that there is an isomorphism of 
$O_{Y,\xi}[I]$-modules 
$$
E_{\xi} \cong \oplus _{1\leq i\leq r} \oplus _{1 \leq j_i \leq f_i} O_{Y,\xi}[I/(^{\gamma_i}H\cap I)] \ .
$$ 
For any $O_{Y,\xi}[I]$-module $M$ and for any integer $1 \leq k \leq e_h/2$, 
let us denote by $M(k)$ the $\chi^k$-component of $M$. 
For each $k$ we deduce  from above the following decomposition
$$ 
E_{\xi}(k) \cong \oplus _{1\leq i\leq r} 
\oplus _{1 \leq j_i \leq f_i} O_{Y,\xi}[I/(^{\gamma_i}H\cap I)](k) \ .
$$ 
Therefore the computation of $d_k(E)$ reduces to the computation of the 
rank $r_{\gamma_i}(k)$ 
over $O_{Y,\xi}$ of  each $O_{Y,\xi}[I/(^{\gamma_i}H\cap I)](k)$. 
It is clear that this module  is 
generated as an  $O_{Y,\xi}$-module by the set 
$\{ \epsilon _{\chi^k}{\bar a}\}$ when ${\bar a}$ runs through 
$I/(^{\gamma_i}H\cap I)$ and where 
$\epsilon _{\chi^k} $ is the idempotent of $O_{Y,\xi}[I]$ 
associated to the character $\chi^k$, namely 
$e^{-1} \sum _{u \in I} \chi^k(u)u^{-1}$.   
For any ${\bar a}$ we have  the equality
$$ 
\epsilon_{\chi^k} {\bar a}=e^{-1}\sum_{u \in I} \chi^k(u)u^{-1}{\bar a} \ .
$$
We now observe that $u^{-1}{\bar a}={\bar a}$ is equivalent to 
$u \in {^{\gamma_i}H}\cap I$. Hence, if for any $1\leq v \leq s$
we denote by $u_v$,  a full set of representatives 
of $I/(^{\gamma_i}H\cap I)$, then  we obtain
$$
\epsilon _{\chi^k}{\bar a}=
e^{-1}\sum_{t \in (^{\gamma_i}H\cap I)}\chi^k(t)
\sum_{1 \leq v \leq s}\chi_h^k(u_v)u_v^{-1}{\bar a} \ .
$$
It follows from the definitions   that the integer $s$ is equal to the ramification index $e_i$. 
It is now clear  from this last equality that $O_{Y,\xi}[I/(^{\gamma_i}H\cap I)](k)$ is different from $0$ if and only if ${\chi^k}$ is trivial when restricted to $ ^{\gamma_i}H\cap I$, namely when $k$ belongs to the set of integers $\{e'_it, 1 \leq t \leq e_i/2 \}$. For any $t$, $1 \leq t \leq e_i/2$, we obtain that
$$
\epsilon _{\chi^{e'_it}}{\bar a}=
e_i^{-1}\sum _{u \in  I/ (^{\gamma_i}H\cap I)}\chi^{e'_it}(u)u {\bar a} \ .
$$
Therefore, for $1 \leq t \leq e_i/2$, we deduce   
that  $O_{Y,\xi}[I/(^{\gamma_i}H\cap I)](e'_it)$ is a free,  rank $1$,  $O_{Y,\xi}$-module with   $\{\epsilon_{\chi^{e'_it}}  \bar 1\}$ as a free basis. Hence  we have proved that $r_{\gamma_i}(k)=r_i(k)=1$ if $e'_i$ divides $k$ and $0$ otherwise.  Therefore we have: 
$$d(E)=\sum _{0 \leq k \leq e/2} kd_k(E))=\sum _{0 \leq k \leq e/2} k\sum _{1  \leq i \leq r} f_i k r_i(k) \ . $$
and thus 
$$
d(E)=\sum _{1  \leq i \leq r}  e'_if_i\sum _{0 \leq k \leq e_i/2} k \ =
\sum _{1  \leq i \leq r} e'_if_i{(e_i^2-1)\over 8} \ .
$$
Since the ramification indices are odd we conclude that in this case 
$$
R(\rho, X)\equiv 
\sum_h (\sum_{\eta_h\rightarrow \xi_h}f(\eta_h){(e(\eta_h)^2-1)\over 8})b_h\ \ \bmod 2\ ,  
$$ where $\eta_h$ ranges over codimension points on $V$ above $\xi_h$. 
The right hand side of this congruence is precisely the divisor obtained in [CNET1] and the divisor obtained by  Serre in [Se2].  
 
\vfill \newpage

\section{Appendix: Hasse-Witt invariants of symmetric complexes}
 
The aim of this note is to indicate how one can define Hasse-Witt 
invariants for symmetric complexes using ideas of  Saito and Walter
(see \cite{TS2} and \cite{Wa} for details). These invariants have been used
in [CNET2] to give a more direct approach to the comparison result
Thm.~0.1 in [CNET1]. They also have been used by Saito
to formulate a conjecture relating invariants of
forms coming from $\ell$-adic and
de Rham cohomology of a projective smooth scheme (see \cite{TS2}).

It should be clear that we do not make any claim to originality. 
The construction of the invariants is due to Walter and Saito.

\subsection{Symmetric complexes}

\noindent We have defined a symmetric bundle $\left(  E,\phi\right)$ 
over a noetherian ${\bf Z}[\frac{1}{2}]$-scheme $Y$ to be 
 a vector bundle $E$ over $Y$ equipped with a
symmetric isomorphism $\phi$ between $E$ and its $Y$-dual $E^{\vee}$, that is
to say $\phi:E\cong E^{\vee}$ and $\phi$ is equal to its transpose $\phi^{\vee}$
after identifying $E$ with $E^{\vee \vee}$. A
{\it symmetric complex} on $Y$ is a symmetric object 
$\left(  P_{\bullet },\phi\right)$ of the derived category 
$\mathcal{D}^{b}\left(  Y\right)  $
of bounded complexes of $Y$-vector bundles. This is a triangulated category
with a dua\-li\-ty which extends $\vee$ : namely the
localisation of the functor which sends a complex $  P_{\bullet}$ 
to the dual complex $P_{\bullet}^{\vee}$. Symmetry is defined using
the natural identification of a complex with its double dual 
(see Sect.~2 of \cite{Ba1}). A symmetric bundle may therefore be viewed as a
symmetric complex concentrated in degree zero. One 
can define a metabolic object in any
triangulated category. In the case of $\mathcal{D}^{b}\left(  Y\right)  $ we
say that $\left(  P_{\bullet},\phi\right)  $ is {\it metabolic} with lagrangian
$L_{\bullet}$ if there is a distinguished triangle
$$
L_{\bullet}\stackrel{i}{\rightarrow}P_{\bullet}\stackrel{i^{\vee}\circ\phi
}{\rightarrow}L_{\bullet}^{\vee}\stackrel{w}{\rightarrow}TL_{\bullet}
$$
in the derived category with the duality condition that $T\left(
w^{\vee}\right)  =w$. One can then define a Witt group for 
$\mathcal{D}^{b}(Y)$, and an object in $\mathcal{D}^{b}(Y)$
is metabolic if and only if it is zero in this Witt group 
(see Thm.~3.5 in \cite{Ba1}).

\begin{prop}
For any symmetric complex $\left(  P_{\bullet},\phi\right)  $ there exists a
sym\-me\-tric bundle $\left(  E^{\prime},\gamma\right)  $ such that the orthogonal
sum $\left(  P_{\bullet},\phi\right)  \perp\left(  E^{\prime},\gamma\right)  $
is metabolic with lagrangian given by $P_{<0}\oplus E^{\prime}.$
\end{prop}

\noindent We briefly sketch a proof of the proposition following indications
of Walter and using the results in \cite{Wa}.
(Note that this reference contains an alternative approach to Balmer's
theorem which identifies the Witt group of $Y$ with the Witt group of
$\mathcal{D}^{b}\left(  Y\right)$(see \cite{Ba2}).)
\medskip

\Proof
The starting point is the existence of a factorisation of
$\phi:P_{\bullet}\rightarrow P_{\bullet}^{\vee}$ as:
\begin{equation}%
\begin{array}
[c]{c}%
P_{\bullet}\\
\downarrow f\\
Q_{\bullet}\\
\downarrow\psi\\
Q_{\bullet}^{\vee}\\
\downarrow f^{\vee}\\
P_{\bullet}^{\vee}
\end{array}
\;\;
\begin{array}
[c]{ccccccccc}%
\cdots & \rightarrow &  P_{1} & \rightarrow &  P_{0} & \stackrel{d_{0}
}{\rightarrow} & P_{-1} & \rightarrow & \cdots\\
&  & \downarrow &  & \downarrow\pi_{1} &  & \| &  & \\
\cdots & \rightarrow &  P_{-1}^{\vee} & \stackrel{\pi_{0}}{\rightarrow} &
E^{\prime \prime} & \rightarrow &  P_{-1} & \rightarrow & \cdots\\
&  & \| &  & \downarrow\lambda &  & \| &  & \\
\cdots & \rightarrow &  P_{-1}^{\vee} & \rightarrow &  E^{\prime \prime 
\vee} & \rightarrow
&  P_{-1} & \rightarrow & \cdots\\
&  & \| &  & \downarrow &  & \downarrow &  & \\
\cdots & \rightarrow &  P_{-1}^{\vee} & \stackrel{d_{0}^{\vee}}{\rightarrow} &
P_{0}^{\vee} & \rightarrow &  P_{1}^{\vee} & \rightarrow & \cdots
\end{array}
\end{equation}
where $\left(  E^{\prime \prime},\lambda \right)  $ is a symmetric bundle
defined using
the acyclicity of the mapping cone of $\phi$. Namely, we factor
$$
P_{-1}^{\vee}\oplus P_{0}\stackrel{\left(
\begin{array}
[c]{cc}
d_{0}^{\vee} & \phi_{0}\\
0 & -d_{0}
\end{array}
\right)  }{\rightarrow}P_{0}^{\vee}\oplus P_{-1}
$$
into a composition of an essential epi $\pi$ with an essential mono $\sigma$
$$
P_{-1}^{\vee}\oplus P_{0}\stackrel{\pi}{\rightarrow}E^{\prime \prime}
\stackrel{\sigma }{\rightarrow}P_{0}^{\vee}\oplus P_{-1}
$$
where $\pi=\left(  \pi_{0},\pi_{1}\right)  ,
\sigma=\left(
\begin{array}[c]{c}
\sigma_{0}\\
\sigma_{1}
\end{array}
\right)  $ (see Lemma 4.2 {\it loc. cit.}). This then defines 
$E^{\prime \prime}$ and
provides a factorisation of $\phi$ as $P_{\bullet}\stackrel{f}{\rightarrow
}Q_{\bullet}
\stackrel{f^{\vee}\circ \psi}{\rightarrow}P_{\bullet}^{\vee}$,
so that there is a
unique $\lambda :E^{\prime \prime}\rightarrow E^{\prime \prime \vee}$ 
with the expected
properties (see Prop. 5.3 {\it loc. cit.}). 
Let $L_{\bullet} : \cdots \rightarrow 0 \rightarrow  0 \rightarrow  P_{-1}
\rightarrow  \cdots$ be the truncation of $P_{\bullet}$ in degrees 
less than zero, then from (5.1) we obtain 
$$
L_{\bullet}\stackrel{u}{\rightarrow} Q_{\bullet} 
\stackrel{\psi}{\rightarrow} Q_{\bullet}^{\vee}
\stackrel{u^{\vee}}{\rightarrow} L_{\bullet}^{\vee} \ .
$$
The complex $L_{\bullet}$ is totally isotropic in $Q_{\bullet},$ i.e.
$u^{\vee}\circ\psi\circ u=0.$ {\it A priori} this holds in 
$\mathcal{D}^{b}(Y)$,  but Prop. 3.1 in [Wa] shows that,
 after possibly changing
$\left(  Q_{\bullet},\psi\right)  $ to $\left(  \widetilde{Q}_{\bullet
},\widetilde{\psi}\right)  $ by an isomorphism in $\mathcal{D}^{b}(Y)$,
 we can lift all maps to maps of chain complexes such that the
resulting chain map lifting $u^{\vee}$ is a split epi. 
If $L_{\bullet}^{\perp}$ denotes the kernel of this map,
then one can consider the subquotient $L_{\bullet}^{\perp}/L_{\bullet},$ on
which $\widetilde{\psi}$ defines a symmetric quasi-isomorphism $q$. In
fact, by Lemma 3.2 of {\it loc. cit.}, $\left[  L_{\bullet}^{\perp}/L_{\bullet
}, q\right]  =\left[  Q_{\bullet},\psi\right]  $ in the derived Witt
group, and hence
$$
\left[  E^{\prime \prime},\lambda\right]  
=\left[  L_{\bullet}^{\perp}/L_{\bullet
}, q\right]  =\left[  Q_{\bullet},\psi\right]  =\left[  P_{\bullet}
,f^{\vee}\psi f\right]  =\left[  P_{\bullet},\phi\right]
$$
and the proposition follows by letting $E^{\prime} = E^{\prime
  \prime}$ and $\gamma = -\lambda$. 
\smallskip

\subsection{Definition of the invariants}

\noindent 
We are now in a position to define the Hasse-Witt invariant of a symmetric
complex. Recall that the total Hasse-Witt invariant of a symmetric bundle
$E=\left( E, \phi\right)  $ over $Y$ is 
$$
w_{t}\left(  E\right)  =%
{\textstyle\sum_{i\geq0}}
w_{i}\left(  E\right)  t^{i}
$$ 
where $w_{i}\left(  E\right)  $ belongs to
$H^{i}\left(  Y\right)  := 
H^{i}\left(  Y_{et},\mathbf{Z/2Z}\right)  $. 
As the main lemma of Sect.~4 shows, 
this invariant does not in general vanish on metabolic bundles.
Recall that if 
 $E$ is a metabolic bundle with lagrangian $V$ of rank $n$ and if
$c_{i}\left(  V\right)  $ is the $i$th Chern class in$\;H^{2i}\left(
Y\right)$, then 

$$
w_{t}\left(  E\right)  =d_{t}\left(  V\right)  
:= {\textstyle\sum_{i=0}^{n}}
\left(  1+\left(  -1\right)  t\right) ^{n-i}c_{i}\left(  V\right)
t^{2i}.
$$
\noindent 
We extend the definition of $d_{t}\left(  -\right)  $ to complexes by
multiplicativity, so that $d_{t}\left(  L_{\bullet}\right)  =%
{\textstyle\prod_{i}}
d_{t}\left(  L_{i}\right)  ^{\left(  -1\right)  ^{i}}.$ Guided by the 
above equality, 
forcing additivity and using the notation of the proposition, we
then recoup Saito's definition of Hasse-Witt classes 
for symmetric complexes by putting 
$$
w_{t}\left(  P_{\bullet}, \phi \right) 
:= w_{t}\left(  E^{\prime}, -\gamma \right)
d_{t}\left(  P_{<0}\right)  \ .
$$
We note that these
Hasse-Witt classes have all the standard properties of characteristic classes.
In particular they are natural with respect to pullback and satisfy the
Whitney sum formula on sums of symmetric complexes; furthermore 
the main lemma extends to metabolic complexes (see \cite{TS2}).

\bigskip

\noindent {\sc Authors' Addresses:}
\medskip

\noindent Ph. C-N. \& B. E.: 
Laboratoire ``Th\'eorie des nombres et algorithmique arithm\'etique'' 
(UMR CNRS 5465), Institut de Math\'ematiques de Bordeaux (FR CNRS 2254),
 Universit\'e Bordeaux~1, 
351, cours de la Lib\'eration, \newline
F-33405 Talence,
France
\smallskip

\noindent e-mail : phcassou@math.u-bordeaux.fr \ \ \
erez@math.u-bordeaux.fr
\bigskip

\noindent M.J. T.: Department of Mathematics, U.M.I.S.T.,
P.O. Box 88, \newline
GB-Manchester M 60 1 QD, England, U.K.
\smallskip

\noindent e-mail : martin.taylor@umist.ac.uk


\vfill 
\newpage

\begin{thebibliography}{C-P-S 123} \addcontentsline{toc}{section}{Bibliography}


\bibitem[Ba1]{Ba1} Balmer P.: {\it  Triangular Witt groups. 
Part I: the 12-term localisation
exact sequences.} K-theory, {\bf 19}(2000), 311-363.

\bibitem[Ba2]{Ba2} Balmer P.: {\it Triangular Witt groups. Part II: 
from usual to derived.} Math. Z. , {\bf 236}(2001), 351-382.

\bibitem[C-E]{CE} Chinburg, T., Erez, B.: {\it
Equivariant Euler-Poincar\'e characteristics and tameness.}
 Journ\'ees Arithm\'etiques, 1991. Ast\'erisque
{\bf 209}(1992), 13, 179--194.


\bibitem[CEPT1]{CEPT-Duke} Chinburg, T., Erez, B., Pappas, G., Taylor, M.J.:
{\it Tame actions of group schemes: integrals and slices.}
 Duke Math. J. {\bf 82}(1996), no. 2, 269--308.

\bibitem[CEPT2]{CEPT} Chinburg, T., Erez, B., Pappas, G., Taylor, M.J.:
{\it  $\epsilon$-constants and the Galois structure of de Rham
cohomology.}
 Ann. of Math. (2){\bf 146}(1997), no. 2, 411--447.

\bibitem[CNET1]{CNET1} Cassou-Nogu\`es, Ph., Erez, B., Taylor, M.J.:
 {\it Invariants of a quadratic form attached to a tame covering of schemes.}
J. Th. des Nombres de Bordeaux  {\bf 12}(2000),  597-660.

\bibitem[CNET2]{CNET2} Cassou-Nogu\`es, Ph., Erez, B., Taylor, M.J.:
{\it Hasse-Witt invariants of symmetric complexes: an example from
geometry.} C. R. Acad. Sci. Paris S\'er. I Math. {\bf 334}(2002),
839-842.

\bibitem[De]{De} Deligne, P.: {\it Les constantes locales de
l'\'equation fonctionnelle de la fonction $L$ d'Artin d'une
repr\'esentation orthogonale.} Invent. Math. {\bf 35}(1976), 299-316.

\bibitem[Dz]{Delz} Delzant, A.: {\it
D\'efinition des classes de Stiefel-Whitney d'un mo\-du\-le
quadratique sur un corps de caract\'eristique diff\'erente de
$2$.} C. R. Acad. Sci. Paris {\bf 255}(1962), 1366--1368.



\bibitem[E]{E2} Erez, B.: {\it 
Geometric trends in Galois module theory.},  116-145, in Galois representations and arithmetic algebraic geometry, eds A.J. Scholl and R.L. Taylor, LMS Lectures Notes {\bf 254}, Cambridge University Press, Cambridge, 1998. 

\bibitem[Es-K-V]{EKV} Esnault, H., Kahn, B., Viehweg, E.: {\it
Coverings with odd ramification and Stiefel-Whitney classes.} J.
reine angew. Math. {\bf 441}(1993), 145--188.

\bibitem[F]{F-SW} Fr\"ohlich, A.: {\it
Orthogonal representations of Galois groups, Stiefel-Whitney
classes and Hasse-Witt invariants.} J. reine angew. Math. {\bf
360}(1985), 84--123.

\bibitem[FQ]{FQ} Fr\"ohlich, A., Queyrut, J.: {\it
On the functional equation of the Artin $L$-function for characters
of real representations.} Invent. Math. {\bf 20}(1973),
125--138.

\bibitem[FT]{FT} Fr\"ohlich, A., Taylor, M.J.: Algebraic number theory. Cambridge studies in advanced Mathematics {\bf 27}. Cambridge University Press , 1991. 

\bibitem[Gl]{Glass} Glass, D.: {\it Invariants associated to orthogonal
$\epsilon$-constants.} Columbia, preprint 2002.

\bibitem[Gr]{Gr-Chern} Grothendieck, A.: {\it
Classes de Chern et repr\'esentations lin\'eaires des groupes
discrets.} Dix Expos\'es sur la Cohomologie des Sch\'emas pp.
215--305, North-Holland, Amsterdam; Masson, Paris, 1968.

\bibitem[Gr-M]{GM} Grothendieck, A., Murre, J.P.:
The tame fundamental group of a formal neighbourhood of a divisor
with normal crossings on a scheme. Lect. Notes in Math., {\bf
208}. Springer-Verlag, Berlin-New York, 1971.

\bibitem[EGA II]{EGAII} Grothendieck, A., Dieudonn\'e, J.:
El\'ements de g\'eom\'etrie al\-g\'e\-bri\-que II. Etude globale 
\'el\'ementaire de quelques classes de morphismes.
Inst. Hautes Etudes Sci. Publ. Math. {\bf 8}(1961).

\bibitem[J1]{J-HW} Jardine, J. F.: {\it
Universal Hasse-Witt classes.}
 Algebraic $K$-theory and algebraic number theory, 83--100,
Contemp. Math., {\bf 83}, Amer. Math. Soc., Providence, RI, 1989.

 \bibitem[J2]{J-HS} Jardine, J. F.: {\it
Higher spinor classes.} Mem. Amer. Math. Soc. {\bf 528}(1994).

\bibitem[K]{K} Kahn, B.: {\it
 Equivariant Stiefel-Whitney classes.}
 J. Pure Appl. Algebra {\bf 97}(1994), no. 2, 163--188.

\bibitem[Kne]{Kne} Knebusch, M.: {\it
Symmetric bilinear forms over algebraic varieties}, Conference on
Quadratic Forms---1976 , pp. 103--283 ed. G. Orzech. Queen's
Papers in Pure and Appl. Math., {\bf 46}, Queen's Univ., Kingston,
Ont., 1977.

\bibitem[Knu]{Knu} Knus, M-A.:
 Quadratic and Hermitian forms over rings.
Grund\-lehren der Math. Wiss. {\bf 294}, Springer-Verlag, Berlin,
1991.




\bibitem[McL]{McL} MacLane, S.: Homology. Grundlehren der Math. Wiss {\bf 114}, Springer-Verlag, Berlin, 1975.

\bibitem[Mi]{Mi} Milne, J.S.:
Etale cohomology. Princeton Mathematical Series, {\bf 33}.
Princeton University Press, Princeton, N.J., 1980.


\bibitem[O]{O2} Ojanguren, M.: The Witt group and the problem
of L\"{u}roth, Dottorato di Ricerca in Mat., Univ. di Pisa,
ETS Editrice, Pisa, 1990.

\bibitem[S1]{TS1} Saito, T.: {\it The sign of the functional equation of the L-function of an orthogonal motive.}
Invent. Math. {\bf 120}(1995), 119-142.

\bibitem[S2]{TS2} Saito, T.: {\it
Note on Stiefel-Whitney class of $\ell$-adic cohomology.}
Preprint, University of Tokyo, 1998.

\bibitem[SGA4]{SGA4} Th\'eorie des topos et cohomologie \'etale des sch\'emas.
S\'eminaire de G\'eom\'etrie Alg\'ebrique du Bois-Marie 1963--1964
(SGA 4). Dirig\'e par M. Artin, A. Grothendieck et J. L. Verdier.
Avec la collaboration de P. Deligne et B. Saint-Donat. Lecture
Notes in Math., {\bf 269, 270, 305}. Springer-Verlag, Berlin-New
York, 1972-73.

\bibitem[Se1]{S-HW} Serre, J-P.: {\it
L'invariant de Witt de la forme ${\rm Tr}(x^2)$.} Com\-ment. Math.
Helv. {\bf 59}(1984), no. 4, 651--676.

\bibitem[Se2]{S-theta} Serre, J-P.: {\it
Rev\^etements \`a ramification impaire et
th\^eta-ca\-ract\'e\-ri\-stiques.} C. R. Acad. Sci. Paris S\'er. I
Math. {\bf 311}(1990), no. 9, 547-552.

\bibitem[Se3]{Serre} Serre, J-P.:
Cohomologie galoisienne. Cinqui\`eme \'edition, r\'evis\'ee et
compl\'et\'ee. Lecture Notes in Mathematics {\bf 5}.
Springer-Verlag, Berlin, 1994.

\bibitem[Se4]{SerreCL} Serre, J-P.:
Corps locaux. Hermann, 1968.

\bibitem[Sn]{Snaith} Snaith, V.P.: {\it
Stiefel-Whitney classes of bilinear forms--a formula of Serre.}
Can. Bull. Math. {\bf 28}(1985), no. 2, 218-222.


\bibitem[Wa]{Wa} Walter, C.: {\it Obstructions to the existence of symmetric
resolutions.} Preprint, Nice, 2001 (see also {\it Grothendieck and Witt
groups of triangulated categories}, Preprint, Nice, 2003).

\end{thebibliography}
\end{document}